\documentclass[12pt]{article}
\usepackage[top=0.95in,bottom=1.5in,left=1.1in,right=1.1in]{geometry}
\usepackage{amsmath,amssymb,amsfonts,epsf,euscript}
\usepackage{graphicx}
\usepackage[colorlinks=true]{hyperref}
\usepackage{float}
\newtheorem{theorem}{Theorem}[section]

\newtheorem{lemma}[theorem]{Lemma}
\newtheorem{corollary}[theorem]{Corollary}
\newtheorem{example}{Example}[section]
\newtheorem{definition}{Definition}[section]
\newtheorem{remark}{Remark}[section]

\renewcommand{\theequation}{\arabic{section}.\arabic{subsection}.\arabic{equation}}
\date{}
\begin{document}
\title{\bf Semi-L\'evy driven continuous-time GARCH process}
\author{M. Mohammadi\thanks{\scriptsize{
Faculty of Mathematics and Computer Science, Amirkabir University of Technology, 
Tehran, Iran. E-mail: m.mohammadiche@aut.ac.ir(M. Mohammadi) and rezakhah@aut.ac.ir(S. 
Rezakhah).}}
\and  S. Rezakhah$^{*}$ \and   N. Modarresi\thanks{\scriptsize{
Department of Mathematics and computer science, Allameh Tabataba'i University, Tehran, 
Iran. E-mail: n.modarresi@atu.ac.ir(N. Modarresi).}}}
\maketitle

\begin{abstract}
In this paper, we study the class of semi-L\'evy driven continuous-time GARCH, denoted by SLD-COGARCH, process. The statistical properties of 
this process are characterized. We show that the state process of such process can be described by a random recurrence equation with the periodic 
random coefficients. We establish sufficient conditions for the existence of a strictly periodically stationary solution of the state process which causes
the volatility process to be strictly periodically stationary. Furthermore, it is shown that the increments with constant length of such SLD-COGARCH
 process are themselves a discrete-time periodically correlated (PC) process. We apply some tests to verify the PC behavior of these increments by 
 the simulation studies. Finally, we show that how well this model fits a set of high-frequency financial data.\\
 
\textit{AMS 2010 Subject Classification}: 60G12, 62M10, 60G51, 60J75, 91B70, 93E15.\\
 
\textit{Keywords}: continuous-time GARCH process; semi-L\'evy process; strictly periodically stationary; periodically correlated.
\end{abstract}

\renewcommand{\theequation}{\arabic{section}.\arabic{equation}}
\section{Introduction}
\setcounter{equation}{0}
The discrete-time ARCH and GARCH processes are often applied to model financial time series \cite{e1, b1}. Kluppelberg et al. \cite{k1} introduced a continuous-time version of the GARCH(1,1), called COGARCH(1,1), process which has the potential to have a better description of changes by determining the underling continuous structure. The examples of such underling structure are L\'evy processes like Brownian motion and Poisson process. They studied L\'evy driven COGARCH(1,1) processes where noises are considered as increments of some L\'evy process and proved the stationarity and also second order properties, under some regularity conditions. Brockwell et al. \cite{b5} generalized the L\'evy driven COGARCH(1,1) process to the higher order. They defined the COGARCH($p,q$) process $(G_{t})_{t\geq0}$ driven by a L\'evy process $(L_{t})_{t\geq0}$ as the solution to the stochastic differential equation
\begin{align*}
dG_{t}=\sqrt{V_{t}}dL_{t},\qquad t>0,\qquad G_{0}=0,
\end{align*}
where the volatility process $(V_{t})_{t\geq0}$ is a continuous-time ARMA($q,p-1$), called CARMA ($q,p-1$), process. They presented the state-space representation of the volatility process as (see e.g. \cite{b4})
\begin{align}\label{equ1.1}
V_{t}=\alpha_{0}+\mathbf{a}^{\prime}{\mathbf{Y}_{t^{-}}},\qquad t>0,\qquad 
V_{0}=\alpha_{0}+\mathbf{a}^{\prime}
{\mathbf{Y}_{0}},
\end{align}
and
\begin{align}\label{equ1.2}
d\mathbf{Y}_{t}=B\mathbf{Y}_{t^{-}}dt+\mathbf{e}(\alpha_{0}+\mathbf{a}^{\prime}\textbf{Y}_{t^{-}})d[L,L]_{t}^{(d)},
\qquad t>0,
\end{align}
where $[L,L]_{t}^{(d)}=\sum_{0<u\leq t}(\Delta L_{u})^{2}$ is the discrete part of the quadratic variation process $([L,L]_{t})_{t\geq0}$ (see Protter \cite[p.66]{p}), $\alpha_{0}>0$ and 
\begin{align*}
B = \begin{pmatrix}
0 & 1 & 0 & \cdots & 0\\
0 & 0 & 1 & \cdots & 0\\
\vdots & \vdots & \vdots & \ddots & \vdots \\
0 & 0 & 0 & \cdots & 1\\
-\beta_{q} & -\beta_{q-1} & -\beta_{q-2} & \cdots & -\beta_{1}\\
\end{pmatrix},\quad
\mathbf{a}=\begin{pmatrix}
\alpha_{1}\\
\alpha_{2}\\
\vdots\\
\alpha_{q}\\
\end{pmatrix},\quad
\mathbf{e}=\begin{pmatrix} 
0\\
0\\
\vdots\\
0\\
1\\
\end{pmatrix},
\end{align*}
with $\alpha_{1},\cdots,\alpha_{p}\in\mathbb{R}, \beta_{1},\cdots, \beta_{q}\in\mathbb{R}, \alpha_{p}\neq0, \beta_{q}\neq0,$  and $\alpha_{p+1}=\cdots=\alpha_{q}=0$. They showed that the state process $(\mathbf{Y}_{t})_{t\geq0}$ can be expressed as a stochastic recurrence equation with the random coefficients. Under some sufficient conditions, they proved that the volatility process $(V_{t})_{t\geq0}$ is strictly stationary which causes the increments with constant length of the process $(G_{t})_{t\geq0}$ to be stationary.

The L\'evy driven COGARCH processes have the restriction that the relevant underling process has the stationary increments. In contrast with such stationary increment processes, processes with periodically stationary increments have a wider application. For example, the accumulated compensation claims from an insurance company in many cases have periodic increments. Specifically, car insurance companies are dealt with more accidental claims at the starting of the new year and in summer than the other parts of the year. As another example, the precipitation in almost all parts of the world have periodic behavior and so the volume of underground waters are higher in some months and lower in some months in each year. The observations of such time series have significant dependency to the previous periods as well. Motivated by these observations, the underling processes with periodically stationary increments, say semi-L\'evy processes, are more prominent than the L\'evy processes. So, we study semi-L\'evy driven COGARCH, denoted by SLD-COGARCH, process in this paper. The semi-L\'evy processes have been extensively studied by Maejima and sato \cite{m1}. Recognizing such semi-L\'evy processes often can be followed by considering some fixed partitions of the corresponding period where inside successive partitions the process has different L\'evy type structure and has similar structure to the corresponding element of the previous periods.

In this paper, we study the statistical properties of the SLD-COGARCH($p,q$) process. Providing a random recurrence equation with periodic random coefficients for the state process $(\mathbf{Y}_{t})_{t\geq0}$ defined by (\ref{equ1.2}), we show that it converges almost surely under some regular conditions. We establish sufficient conditions for the existence of a strictly periodically stationary solution of the state process $(\mathbf{Y}_{t})_{t\geq0}$ which causes the volatility process $(V_{t})_{t\geq0}$ defined by (\ref{equ1.1}) to be strictly periodically stationary. Furthermore, the increments with constant length of the SLD-COGARCH($p,q$) process are a discrete-time periodically correlated (PC) process. We use some tests to verify the PC behavior of the increments. Finally, we apply the introduced model to some real data set.

This paper is organized as follows. In section 2, we introduce the SLD-COGARCH processes. We also present the structure of the semi-L\'evy process and obtain its characteristic function. Section 3 is devoted to some sufficient conditions that cause the volatility process to be strictly periodically stationary. We also show that the increments with constant length of such SLD-COGARCH process form a PC process. In section 4, we illustrate the results with simulations. In section 5, we apply the model to a set of high-frequency financial data and show that this process has the consistent properties. All proofs are presented in Section 6.

\section{Semi-L\'evy driven COGARCH process}
\setcounter{equation}{0}
In this section, we describe the structure of semi-L\'evy compound Poisson process in subsection 2.1. Following the method of Brockwell et al. \cite{b5}, we present the structure of corresponding semi-L\'evy compound Poisson process as the underling process of the COGARCH process in subsection 2.2.
\subsection{Semi-L\'evy process}
We present the definition of semi-L\'evy Poisson and semi-L\'evy compound Poisson processes. Then, we drive the characteristic function of the semi-L\'evy compound Poisson process. 
\begin{definition} \textnormal{\textbf{(Semi-L\'evy Poisson process)}}\\
Let $A_{j}=(s_{_{j-1}},s_{_{j}}]$  where $0=s_{_{0}}<s_{_{1}}<\cdots$ constitutes a partition of positive real line 
and $|A_{j}|=|A_{j+d}|=s_{_{j}}-s_{_{j-1}}$
 for some integer $d$ 
and all $j\in\mathbb{N}$
 and also $\tau=\sum_{j=1}^{d}|A_{j}|$.
 Then, a non-homogeneous Poisson process $\big(N(t)\big)_{t\geq0}$
 is called a semi-L\'evy Poisson process with period $\tau>0$ and intensity function $\lambda(u)$ 
if $E\big(N(t)\big)=\Lambda(t)$
where 
\begin{align}\label{equ2.1}
\Lambda(t)=\int_{0}^{t}\lambda(u)du
\end{align}
and $\lambda(t)=\sum_{j=1}^{d}\lambda_{j}I_{\mathfrak{D}_{j}}(t)$ for $t,\lambda_{j}\geq0$ and $\mathfrak{D}_{j}=\bigcup_{k=0}^{\infty}A_{j+kd}$.
\end{definition}
One can easily verify that a semi-L\'evy Poisson process has periodically stationary increments.
\begin{definition} \textnormal{\textbf{(Semi-L\'evy compound Poisson process)}}\\
Let $\big(N(t)\big)_{t\geq0}$ be a semi-L\'evy Poisson process with some period $\tau>0$. Then, the semi-L\'evy compound Poisson, denoted by SLCP, process $(S_{t})_{t\geq0}$ is defined as 
\begin{align}\label{equ2.2}
S_{t}=\delta t+\sum_{n=1}^{N(t)}Z_{n}
\end{align}
where $\delta\in\mathbb{R}$ and $Z_{n}=\sum_{j=1}^{d}Z_{n}^{(j)}I_{\lbrace \Upsilon_{_{n}}\in\mathfrak{D}_{j}\rbrace}$ in which $\Upsilon_{_{n}}$ is the arrival time of $n^{th}$ jump $Z_{n}$, $\mathfrak{D}_{j}=\bigcup_{k=0}^{\infty}A_{j+kd}$ and $Z_{n}^{(j)}$ are independent and identically distributed (i.i.d.) random variables with distribution $F_{j}, j=1, \cdots, d$, such that $\int_{\mathbb{R}}z^{2}F_{j}(dz)<\infty$.
\end{definition}
The characteristic function of this process is derived from the following Lemma.
\begin{lemma}
Let $(S_{t})_{t\geq0}$ be the SLCP process introduced by Definition 2.2. Then its characteristic function has the following L\'evy-Khintchine representation
\begin{align}\label{equ2.3}
E[e^{iuS_{t}}]=exp\Big[iu\gamma(t)+\int_{\mathbb{R}}(e^{iu\mathit{z}}-1-iu\mathit{z} 
I_{\lbrace|z|\leq1\rbrace})\nu_{t}(dz)\Big],
\end{align}
where 
\begin{align*}
\gamma(t)=\delta t+\int_{|z|\leq1}z\nu_{t}(dz)
\end{align*}
and
\begin{align*}
\nu_{t}(dz)=\sum_{j=1}^{r-1}(m+1)\lambda_{j}|A_{j}|F_{j}(dz)+\sum_{j=r}^{d}m\lambda_{j}|A_{j}|F_{j}(dz)+
\lambda_{r}(t-s_{_{r-1+md}})F_{r}(dz),
\end{align*}
in which $m=[\frac{t}{\tau}]$ and $r\in\lbrace1,\cdots, d\rbrace$.
\end{lemma}
Proof: see Appendix, P1.
\begin{remark}
(i) By the fact that the numbers of jumps of the SLCP process $(S_{t})_{t\geq0}$ in the time interval $[0,t]$ are finite and the size of each jump is finite as well, we deduce that the process $(S_{t})_{t\geq0}$ is a finite variation process, see \cite[p.14]{r1}. Hence, the process $(S_{t})_{t\geq0}$ is a semimartingale \cite[p.55]{p}.

(ii) By Definition 2.2, we have that if there is a jump at point $u$, $0\leq u\leq t$, then 
\begin{align}\label{equ2.4}
\Delta S_{u}:=S_{u}-S_{u^{-}}&=\big[\delta u+\sum_{n=1}^{N(u)}Z_{n}\big]-\big[\delta 
u^{-}+\sum_{n=1}^{N(u^{-})}Z_{n}\big]=Z_{N(u)},
\end{align}
otherwise $\Delta S_{u}=0$.  Thus, by the Theorem 22 of Protter \cite[p.66]{p}, if there is a jump at point $u$ then $\Delta[S,S]_{u}=|Z_{N(u)}|^{2}$ and otherwise $\Delta[S,S]_{u}=0$. So, the quadratic variation of the process $(S_{t})_{t\geq0}$ is given by
\begin{align}\label{equ2.5}
[S,S]_{t}=\sum_{0\leq u\leq t}\Delta[S,S]_{u}=\sum_{0\leq u\leq t}|Z_{N(u)}|^{2}=\sum_{n=1}^{N(t)}
|Z_{n}|^{2}.
\end{align}
\end{remark}

\subsection{Semi-L\'evy driven COGARCH process}
The representation of the L\'evy driven COGARCH($p,q$) process is based on the structure of the GARCH($p,q$) process $(\xi_{n})_{n\in\mathbb{N}^{0}}$ as
\begin{align*}
\xi_{n}&=\sqrt{V_{n}}\varepsilon_{n}\\
V_{n}&=\alpha_{0}+\alpha_{1}\xi_{n-1}^{2}+\cdots+\alpha_{p}\xi_{n-
p}^{2}+\beta_{1}V_{n-1}+\cdots+\beta_{q}V_{n-q},\quad n\geq s,
\end{align*}
where $s=max\lbrace p, q\rbrace$ and $(\varepsilon_{n})_{n\in\mathbb{N}^{0}}$ is a sequence of i.i.d. random variables. The process $(\xi_{n})_{n\in\mathbb{N}^{0}}$ can be interpreted as the log return of the asset price process $(G_{n})_{n\in\mathbb{N}^{0}}$, say $ln(G_{n+1}/G_{n})$. The volatility process $(V_{n})_{n\in\mathbb{N}^{0}}$ is as an ARMA($q,p-1$) process driven by the noise sequence $(V_{n-1}\varepsilon_{n-1}^{2})_{n\in\mathbb{N}}$. Motivated by this, the dynamic of the L\'evy driven COGARCH($p,q$) process for the logarithm of the asset price process $(G_{t})_{t\geq0}$ is defined as
\begin{align*}
dG_{t}=\sqrt{V_{t}}dL_{t}, \qquad t>0, \qquad G_{0}=0,
\end{align*}
where $(L_{t})_{t\geq0}$ is a L\'evy process and $(V_{t})_{t\geq0}$ is a continuous-time analog of the self-exciting ARMA($q,p-1$) process $(V_{n})_{n\in\mathbb{N}^{0}}$, see Brockwell et al. \cite{b5}.

We generalize the L\'evy driven COGARCH($p,q$) process by replacing the L\'evy process $(L_{t})_{t\geq0}$ with the SLCP process $(S_{t})_{t\geq0}$.

\begin{definition}\textnormal{\textbf{\big(SLD-COGARCH($p,q$)\big)}}\\
Let $(S_{t})_{t\geq0}$ be the SLCP process with some period $\tau>0$ defined by (\ref{equ2.2}). Then the left-continuous volatility process $(V_{t})_{t\geq0}$ is defined as 
\begin{align}\label{equ2.6}
V_{t}=\alpha_{0}+\mathbf{a}^{\prime}{\mathbf{Y}_{t^{-}}},\qquad t>0,\qquad 
V_{0}=\alpha_{0}+\mathbf{a}^{\prime}
{\mathbf{Y}_{0}},
\end{align}
where the state process $(\mathbf{Y}_{t})_{t\geq0}$ is the unique c\`adl\`ag solution of the stochastic differential equation
\begin{align}\label{equ2.7}
d\mathbf{Y}_{t}=B\mathbf{Y}_{t^{-}}dt+\mathbf{e}(\alpha_{0}+\mathbf{a}^{\prime}\textbf{Y}_{t^{-}})d[S,S]_{t},
\qquad t>0,
\end{align}
with an initial value $\mathbf{Y}_{0}$ which is independent of the driving process $(S_{t})_{t\geq0}$, in which $[S,S]_{t}$ is the quadratic variation of the process $(S_{t})_{t\geq0}$ and the ($q\times q$)-matrix $B$ and vectors $\mathbf{a}$ and $\mathbf{e}$ are defined as
\begin{align}\label{equ2.8}
B = \begin{pmatrix}
0 & 1 & 0 & \cdots & 0\\
0 & 0 & 1 & \cdots & 0\\
\vdots & \vdots & \vdots & \ddots & \vdots \\
0 & 0 & 0 & \cdots & 1\\
-\beta_{q} & -\beta_{q-1} & -\beta_{q-2} & \cdots & -\beta_{1}\\
\end{pmatrix},\quad
\mathbf{a}=\begin{pmatrix}
\alpha_{1}\\
\alpha_{2}\\
\vdots\\
\alpha_{q}\\
\end{pmatrix},\quad
\mathbf{e}=\begin{pmatrix} 
0\\
0\\
\vdots\\
0\\
1\\
\end{pmatrix},
\end{align}
where $p$ and $q$ are integers such that $q\geq p\geq1, \alpha_{0}>0, \alpha_{1},\cdots,\alpha_{p}\in\mathbb{R}, \beta_{1},\cdots,\beta_{q}\in\mathbb{R}, \alpha_{p}\neq0, \beta_{q}\neq0,$ and $\alpha_{p+1}=\cdots=\alpha_{q}=0$.\\
If the volatility process $(V_{t})_{t\geq0}$ defined by (\ref{equ2.6}) is almost surely non-negative, then $(G_{t})_{t\geq0}$ given by
\begin{align}\label{equ2.9}
dG_{t}=\sqrt{V_{t}}dS_{t},\qquad t>0,\quad G_{0}=0,
\end{align}
is called the SLCP driven COGARCH(p,q), denoted by SLD-COGARCH(p,q), process with parameters $B$, $\mathbf{a}$, $\alpha_{0}$ and the driving process $(S_{t})_{t\geq0}$.
\end{definition}

We present the definition of matrix norms and some notation that will be used throughout the paper.
\begin{definition}
Let $||\cdot||_{r}$ denote the vector $L^{r}$-norm in $\mathbb{C}^{q}$, for $r\in[1,\infty]$. Then the matrix $L^{r}$-norm of the $(q\times q)-$matrix $C$ is defined as
\begin{align}\label{equ2.10}
||C||_{r}=\sup_{\mathbf{c}\in\mathbb{C}^{q}\setminus\lbrace0\rbrace}\frac{||C\mathbf{c}||_{r}}
{||\mathbf{c}||_{r}}.
\end{align}
\end{definition}
The matrix $L^{r}$-norms of the $L^{1}$, $L^{2}$ and $L^{\infty}$  are called the maximum absolute column sum norm, spectral norm and maximum absolute row sum norm, respectively. 
\begin{definition}
Let $P$ be a matrix such that $P^{-1}BP$ is a diagonal matrix, where B is defined by (\ref{equ2.8}). Then the natural matrix norm of  $C$ is
\begin{align*}
||C||_{P,r}=||P^{-1}CP|| _{r}.
\end{align*}
\end{definition}
The natural matrix norm $||\cdot||_{P,r}$  corresponding to the natural vector norm is defined as
\begin{align}\label{equ2.11}
||\mathbf{c}||_{P,r}=||P^{-1}\mathbf{c}|| _{r}, \qquad \mathbf{c}\in\mathbb{C}^{q}.
\end{align}
The eigenvalues of the matrix B are denoted by $\eta_{1}, \cdots, \eta_{q}$ and $\eta:=\max_{i=1,\cdots,q}\mathcal{R}e(\eta_{i})$, where $\mathcal{R}e(\eta_{i})$ is the real part of the eigenvalue $\eta_{i}$. The ($q\times q$)-identity matrix is denoted by $I_{q}$ or simply $I$.
\section{Periodic stationarity of the increments}
\setcounter{equation}{0}
In this section, we show that under some conditions, the volatility process $(V_{t})_{t\geq0}$ defined by (\ref{equ2.6}) is strictly periodically stationary with period $\tau$. The increments with constant length of the SLD-COGARCH process are a periodically correlated (PC) process. Furthermore, under a necessary and sufficient condition, the volatility process $(V_{t})_{t\geq0}$ is non-negative.

The following theorem shows that the state process of the SLD-COGARCH process defined by (\ref{equ2.7}) satisfies a multivariate random recurrence equation with periodic random coefficients.
\begin{theorem}
Let $(\mathbf{Y}_{t})_{t\geq0}$ be the state process of the SLD-COGARCH process. Then the following results are hold.
\begin{itemize}
\item[(a)] There exists a random $(q\times q)-$matrix $J_{s,t}$ and a random vector $\mathbf{K}_{s,t}$ such that
\begin{align}\label{equ3.1}
\mathbf{Y}_{t}=J_{s,t}\mathbf{Y}_{s}+\mathbf{K}_{s,t},\qquad 0\leq s\leq t.
\end{align}
\item[(b)] $\big(J_{s,t},\mathbf{K}_{s,t}\big)$ are periodic in both indices with the same period $\tau$. So, the distribution of $\big(J_{s,t}, \mathbf{K}_{s,t}\big)$ is invariant under time changes of $s,t$ to $s+\tau,t+\tau$.
\item[(c)] $\big(J_{s,t},\mathbf{K}_{s,t}\big)$ and $\big(J_{r,u}, \mathbf{K}_{r,u}\big)$ are independent for $0\leq s\leq t\leq r\leq u$. $\mathbf{Y}_{s}$ is also independent of $\big(J_{s,t},\mathbf{K}_{s,t}\big)$.
\end{itemize}
\end{theorem}
Proof: see Appendix, P2.

In the next theorem, we consider some mild conditions to the state process that converges in distribution to a finite random vector.
\begin{theorem}
Let $(\mathbf{Y}_{t})_{t\geq0}$ be the state process of the SLD-COGARCH process that is satisfied in the followings
\begin{itemize}
\item[(i)] all the eigenvalues of the matrix B are distinct and have strictly negative real parts,
\item[(ii)] there exists some $r\in[1,\infty]$ such that for $t\in[0,\tau]$ and $j=1, \cdots, d$
\begin{align}\label{equ3.2}
\int_{\mathbb{R}}log\big(1+||P^{-1}\mathbf{e}\mathbf{a}^{\prime}P||_{r}z^{2}\big)F_{j}(dz)
<-\frac{\eta\tau}{\Lambda(t+\tau)-\Lambda(t)},
\end{align}
where $P$ is a matrix which causes $P^{-1}BP$ is diagonal.
\end{itemize}
Then, for the fixed $t\in[0,\tau)$, $\mathbf{Y}_{t+m\tau}$ converges in distribution to some vector $\mathbf{U}^{(t)}$, as $m\rightarrow\infty$. Furthermore, $\mathbf{U}^{(t)}$ has a unique distribution that is satisfied in equation
\begin{align}\label{equ3.3}
\mathbf{U}^{(t)}\overset{d}{=}J_{t,t+\tau}\mathbf{U}^{(t)}+\mathbf{K}_{t,t+\tau},
\end{align}
where $\overset{d}{=}$ denotes the equality in 
distribution and $\mathbf{U}^{(t)}$ is independent of $\big(J_{t,t+\tau},\mathbf{K}_{t,t+\tau}\big)$. 
\end{theorem}
Proof: see Appendix, P3.
\begin{corollary}
Under the conditions (i) and (ii) of Theorem 3.2, the following results are hold.
\begin{itemize}
\item[(a)] If $\mathbf{Y}_{t}\overset{d}{=}\mathbf{U}^{(t)}$ for any $t\in[0, \tau)$, then $\mathbf{Y}_{t}$ and $V_{t}$ defined by (\ref{equ2.6}) are strictly periodically stationary with the period $\tau$. In other words, for any $n\in\mathbb{N}$ and $t_{_{1}}, t_{_{2}}, \cdots, t_{_{n}}\geq0$
\begin{align*}
\big(\mathbf{Y}_{t_{_{1}}}, \mathbf{Y}_{t_{_{2}}}, \cdots, \mathbf{Y}_{t_{_{n}}}\big)\overset{d}{=}
\big(\mathbf{Y}_{t_{_{1}}+\tau}, \mathbf{Y}_{t_{_{2}}+\tau}, \cdots, \mathbf{Y}_{t_{_{n}}+\tau}\big)
\end{align*}
and
\begin{align*}
\big(V_{t_{_{1}}}, V_{t_{_{2}}}, \cdots, V_{t_{_{n}}}\big)\overset{d}{=}\big(V_{t_{_{1}}+\tau}, V_{t_{_{2}}+\tau}, \cdots, 
V_{t_{_{n}}+\tau}\big).
\end{align*}
\item[(b)] Increments with the constant length of the SLD-COGARCH process $(G_{t})_{t\geq0}$ is a PC process. In other words, $G_{t}^{(l)}:=\int_{t}^{t+l}\sqrt{V_{u}}dS_{u}$ for any $t, h\geq0$ and fixed $l>0$ satisfies
\begin{align*}
E\big(G_{t}^{(l)}\big)=E\big(G_{t+\tau}^{(l)}\big)
\end{align*}
and
\begin{align*}
cov\big(G_{t}^{(l)},G_{t+h}^{(l)}\big)=cov\big(G_{t+\tau}^{(l)},G_{t+h+\tau}^{(l)}\big).
\end{align*}
\end{itemize}
\end{corollary}
Proof: see Appendix, P4.

One of the striking features of high-frequency financial data is that the returns have negligible correlation while the squared returns are significantly
correlated. These data typically show a PC structure in their squared logarithm returns (known as squared log return). 
The next corollary shows that under some conditions, the increments with constant length of the SLD-COGARCH 
which were defined in Corollary 3.3 (b) can capture these features.
\begin{corollary}
Let the conditions of Theorem 3.2 hold and the driving process $(S_{t})_{t\geq0}$ defined by (2.2) has zero mean. Then increments with constant length $G_{t}^{(l)}$ for $t\geq0$ and $h\geq l>0$ satisfies 
\begin{align*}
&(a)\quad E\big(G_{t}^{(l)}\big)=0, \qquad cov\big(G_{t}^{(l)},G_{t+h}^{(l)}\big)=0,\\ \\
&(b)\quad E\Big((G_{t}^{(l)})^{2}\Big)=E\Big((G_{t+\tau}^{(l)})^{2}\Big),\qquad 
cov\Big((G_{t}^{(l)})^{2},(G_{t+h}^{(l)})^{2}\Big)=cov\Big((G_{t+\tau}^{(l)})^{2},(G_{t+h+\tau}^{(l)})^{2}\Big).
\end{align*}
\end{corollary}
Proof: see Appendix, P5.

The following theorem represents the necessary and sufficient conditions for non-negativity of the volatility process $V_{t}$ for any $t\geq0$. 
\begin{theorem}
Let $(\mathbf{Y}_{t})_{t\geq0}$ be the state process of the SLD-COGARCH process with parameters $B$, $\mathbf{a}$ and $\alpha_{0}>0$. Suppose that $\gamma\geq-\alpha_{0}$ be a real constant that satisfies in 
\begin{align}\label{equ3.4}
\mathbf{a}^{\prime}e^{Bt}\mathbf{e}\geq0\quad \forall t\geq0,
\end{align}
and
\begin{align}\label{equ3.5}
\mathbf{a}^{\prime}e^{Bt}\mathbf{Y}_{0}\geq\gamma \quad a.s.\quad \forall t\geq0.
\end{align}
Then, with probability one, $V_{t}\geq\alpha_{0}+\gamma\geq0$ for $t\geq0$.\\
Conversely, if either (\ref{equ3.5}) fails or (\ref{equ3.5}) holds with $\gamma>-\alpha_{0}$ and (\ref{equ3.4}) fails, then there exists a SLCP process $(S_{t})_{t\geq0}$ and $t\geq0$ such that $P(V_{t} < 0)>0$.
\end{theorem}
Proof: see Appendix, P6.
\section{Simulation}
\setcounter{equation}{0}
In this section, we verify the theoretical results concerning the PC structure of the increments of the SLD-COGARCH process $(G_{t})_{t\geq0}$ is defined in (\ref{equ2.9}) through the simulation. For this, we first simulate the driving semi-L\'evy process $(S_{t})_{t\geq0}$ with period $\tau$ introduced by Definition 2.2. Second, we present the discretized version of the state process $(\mathbf{Y}_{t})_{t\geq0}$ defined by (\ref{equ2.7}) at the time of jump points. We generate the volatility process $(V_{t})_{t\geq0}$ defined by (\ref{equ2.6}) and the process $(G_{t})_{t\geq0}$ at the time of jump points. Then, we evaluate a sampled process from the simulated values of the process $(G_{t})_{t\geq0}$. Finally, we apply some tests to demonstrate the PC structure of the increments for the sampled process.

For simulating the driving process $(S_{t})_{t\geq0}$ with the underlying semi-L\'evy Poisson process $\big(N(t)\big)_{t\geq0}$, we consider $\Upsilon_{_{j}}$ as the time of $j^{th}$ jump $\big(N(t)\big)_{t\geq0}$ and $\Upsilon_{_{0}}=0$. First, by Definition 2.1, we generate the arrival times $\Upsilon_{_{1}}, \Upsilon_{_{2}}, \cdots$ of the process $\big(N(t)\big)_{t\geq0}$ by the following algorithm.
\begin{enumerate}
\item 
Consider the positive value $\tau$ as the length of the period of the process $(N_{t})_{t\geq0}$ and some integer $m$ as the number of period intervals for the simulation.
\item
Choose the positive integer $d$ as the number of elements of partition in each period interval. 
\item
Consider positive real numbers $l_{1}, \cdots, l_{d}$ such that $\tau=\sum_{j=1}^{d}l_{j}$ and a partition of first period interval $(0,\tau]$ by $A_{1}, A_{2}, \cdots, A_{d}$ where $A_{j}=(s_{_{j-1}},s_{_{j}}]$, $s_{_{0}}=0$ and $s_{_{j}}=\sum_{i=1}^{j}l_{i}$, $j=1, 2, \cdots, d$. Elements of partition of $k^{th}$ period interval $\big((k-1)\tau,k\tau\big]$ are $A_{1+(k-1)d}, \cdots, A_{kd}$ where $|A_{j+(k-1)d}|=|A_{j}|=l_{j}$ for $k=2, \cdots, m$ and $j=1, \cdots, d$.
\item
Let the positive real numbers $\lambda_{1}, \cdots, \lambda_{d}$ be occurrence rates corresponding to the increments of $(N(t))_{t\geq0}$ on $A_{1}, \cdots, A_{d}$. It is also followed from Definition 2.1 that the process $(N(t))_{t\geq0}$ has the occurrence rate $\lambda_{j}$ on $A_{j+(k-1)d}$ for $k=2, \cdots, m$ and $j=1, \cdots, d$.
\item
Generate an independent sequence of Poisson random variables $P_{i}$ with parameter $\lambda_{i}l_{i}$ on $A_{i}$ for $i=1, \cdots, md$ as $p_{1}, \cdots, p_{md}$.
\item
Generate $p_{i}$ independent samples from uniform distribution $U(s_{_{i-1}},s_{_{i}}]$ for $i=1, \cdots,$ $md$. Then, sort these samples and denote these ordered samples by $\Upsilon_{_{1}}, \cdots, \Upsilon_{_{md}}$.
\item
Generate the successive jump size $Z_{n}$ independently and with distribution $F_{j}$ if corresponding arrival time belongs to $\mathfrak{D}_{j}=\bigcup_{k=0}^{\infty}A_{j+kd}$  for $ j=1, \cdots, d$ and evaluate $S_{t}$ by (\ref{equ2.2}).
\end{enumerate}

Now, by the following steps, we simulate the SLD-COGARH($p,q$) process introduced by Definition 2.3.
\begin{enumerate}
\item
Set integer valued $p$ and $q$ such that $q\geq p\geq1.$
\item
Choose the real parameters $\beta_{1}, \cdots, \beta_{q}$ and $\alpha_{1}, \cdots, \alpha_{p}$ and $\alpha_{0}>0$ such that $\alpha_{p}\neq0$, $\beta_{q}\neq0$, $\alpha_{p+1}=\cdots=\alpha_{q}=0$. Furthermore, the eigenvalues of the matrix $B$, $\eta_{1},\cdots,\eta_{q}$, are distinct and have strictly negative real parts and conditions (\ref{equ3.2}), (\ref{equ3.4}) and (\ref{equ3.5}) are satisfied.\\
To verify the condition (3.2), define the matrix $P$ as follows
\begin{align*}
P = \begin{pmatrix}
1 & \cdots & 1\\
\eta_{1} & \cdots & \eta_{q}\\
\vdots & \cdots & \vdots \\
\eta_{1}^{q-1} & \cdots & \eta_{q}^{q-1}\\
\end{pmatrix}.
\end{align*}
\item
After generating the arrival times $\Upsilon_{_{n}}, n=1, \cdots, md$, by the above algorithm, generate the state process $(\mathbf{Y}_{\Upsilon_{_{n}}})_{n\in\mathbb{N}}$ in (\ref{Th3.1(5.5)}) with an initial value $\mathbf{Y}_{\Upsilon_{0}}$
\begin{align*}
\mathbf{Y}_{\Upsilon_{_{n}}}&=(I+Z_{n}^{2}\mathbf{e}\mathbf{a}^{\prime})e^{B(\Upsilon_{_{n}}-\Upsilon_{_{n-1}})}
\mathbf{Y}_{\Upsilon_{_{n-1}}}+\alpha_{0}Z_{n}^{2}\mathbf{e}\\
&=e^{B\big(\Upsilon_{_{n}}-\Upsilon_{_{n-1}}\big)}\mathbf{Y}_{\Upsilon_{_{n-1}}}+\mathbf{e}\Big(\alpha_{0}+
\mathbf{a}^{\prime}e^{B\big(\Upsilon_{_{n}}-\Upsilon_{_{n-1}}\big)}\mathbf{Y}_{\Upsilon_{_{n-1}}}\Big)Z^{2}_{n}.
\end{align*}
\item
Similar to (\ref{Th3.1(5.2)}) $\mathbf{Y}_{t}$ for $t\in[\Upsilon_{_{n-1}},\Upsilon_{_{n}}^{-}]$ can be written as $\mathbf{Y}_{t}=e^{B(t-\Upsilon_{n-1})}\mathbf{Y}_{\Upsilon_{n-1}}$. 
So, by (\ref{equ2.6}), generate the process $(V_{\Upsilon_{n}})_{n\in\mathbb{N}^{0}}$ with 
\begin{align}\nonumber
V_{\Upsilon_{_{n}}}&=\alpha_{0}+\mathbf{a}^{\prime}\mathbf{Y}_{\Upsilon_{_{n}}^{-}}\\ \label{equ4.2}
&=\alpha_{0}+\mathbf{a}^{\prime}e^{B(\Upsilon_{_{n}}-\Upsilon_{_{n-1}})}\mathbf{Y}_{\Upsilon_{_{n-1}}}.
\end{align}
\item
Since the driving process $(S_{t})_{t\geq0}$ has one jump at time $\Upsilon_{n}$ over $[\Upsilon_{n-1},\Upsilon_{n}]$, it follows from (\ref{equ2.9}) that 
\begin{align}\nonumber
G_{\Upsilon_{n}}-G_{\Upsilon_{n-1}}=\int_{0}^{\Upsilon_{n}}\sqrt{V_{u}}dS_{u}-\int_{0}^{\Upsilon_{n-1}}\sqrt{V_{u}}dS_{u}
&=\int_{\Upsilon_{n-1}}^{\Upsilon_{n}}\sqrt{V_{u}}dS_{u}\\ \label{equ4.3}
&=\sqrt{V_{\Upsilon_{n}}}Z_{n}.
\end{align}
Then generate the process $(G_{\Upsilon_{n}})_{n\in\mathbb{N}^{0}}$ by (\ref{equ4.3}) such that $G_{\Upsilon_{0}}=G_{0}=0$.
\item
Finally, using the values of $V_{\Upsilon_{n}}$ and $G_{\Upsilon_{n}}$ provided by the previous steps, evaluate the sampled processes $(V_{il})_{i\in\mathbb{N}^{0}}$ and $(G_{il})_{i\in\mathbb{N}^{0}}$ for $l=\tau/\varrho$ where $\varrho$ is some integer by the following steps:
\begin{itemize}
\item
For $il\in[\Upsilon_{_{n-1}},\Upsilon_{_{n}}), n\in\mathbb{N}$, it follows from (\ref{Th3.1(5.2)}) that $\mathbf{Y}_{il}=e^{B(il-\Upsilon_{_{n-1}})}\mathbf{Y}_{\Upsilon_{_{n-1}}}$. So, using (\ref{equ2.6}), we have
\begin{align*}
V_{il}&=\alpha_{0}+\mathbf{a}^{\prime}\mathbf{Y}_{il^{-}}\\
&=\alpha_{0}+\mathbf{a}^{\prime}e^{B(il-\Upsilon_{n-1})}\mathbf{Y}_{\Upsilon_{n-1}}.
\end{align*}
Note that if $il=\Upsilon_{_{n-1}}$, by (\ref{equ4.2})
\begin{align*}
V_{il}=V_{\Upsilon_{_{n-1}}}=\alpha_{0}+\mathbf{a}^{\prime}e^{B(\Upsilon_{_{n-1}}-\Upsilon_{_{n-2}})}
\mathbf{Y}_{\Upsilon_{_{n-2}}}.
\end{align*}
\item
Using the fact that the driving process $(S_{t})_{t\geq0}$ has no jump over $[\Upsilon_{_{n-1}},\Upsilon_{_{n}})$, for $n\in\mathbb{N}$, it follows from (\ref{equ2.9}) for $il\in[\Upsilon_{_{n-1}},\Upsilon_{_{n}})$
\begin{align}\label{equ4.4}
G_{il}-G_{\Upsilon_{_{n-1}}}=\int_{\Upsilon_{_{n-1}}}^{il}\sqrt{V_{u}}dS_{u}=0.
\end{align}
Thus, $G_{il}=G_{\Upsilon_{n-1}}$. If $il=\Upsilon_{_{n-1}}$, then it follows from (\ref{equ4.3}) that
\begin{align}\label{equ4.5}
G_{il}=G_{\Upsilon_{_{n-1}}}=G_{\Upsilon_{_{n-2}}}+\sqrt{V_{\Upsilon_{_{n-1}}}}Z_{n-1}.
\end{align}
\end{itemize}
\end{enumerate}

To detect the PC structure of a process, it is shown that the proposed sample spectral coherence statistic can 
be used to test whether a discrete-time process is PC \cite{h1, h2}.  
In this method, for $n$ samples $X_{0}, X_{1}, \cdots, X_{n-1}$, and a fixed $M$, the following sample spectral coherence statistic is plotted
\begin{align*}
|\hat{\gamma}(P, Q ,M)|^{2}=\frac{|\sum_{m=0}^{M-1}d_{X}(\zeta_{_{P+m}})\overline{d_{X}
(\zeta_{_{Q+m})}}|^{2}}{\sum_{m=0}^{M-1}|d_{X}(\zeta_{_{P+m}})|^{2}\sum_{m=0}^{M-1}|d_{X}(\zeta_{_{Q+m}})|^{2}},
\end{align*}
where $P, Q=0, \cdots, n-1$, $d_{X}(\zeta_{_{P}})=\sum_{k=0}^{n-1}X_{k}e^{ik\zeta_{_{P}}}$ and $\zeta_{_{P}}=\frac{2\pi P}{n}$ for $P=0,\cdots, n-1$.\\
The perception of this sample spectral coherence is aided by plotting the coherence values only at points where the following $\alpha-$threshold is exceeded \cite[p.310]{h2}
\begin{align*}
x_{\alpha}=1-e^{log(\alpha)/(M-1)}.
\end{align*}
For a PC time series with period $\varrho$, it is expected that the sample spectral coherence statistic has a significant value on pairs $(P,Q)$ for $|P-Q|=kn/\varrho, k=0, \cdots, \varrho-1$. The plot of the support of the sample spectral coherence statistic can be used to identify the type of model and analysis of the time series $X_{0}, X_{1}, \cdots, X_{n-1}$, see \cite{d2}:
\begin{itemize}
\item
If only the main diagonal appears, then the time series $X_{i}$ is a stationary time series.
\item
If there are some significant values of statistic and they seem to lie along the parallel equally spaced diagonal lines, then 
the time series $X_{i}$ is PC with period $\varrho=n/L$, where $L$ is the line spacing.
\item
If there are some significant values of statistic but they occur in some non-regular places, then 
the time series $X_{i}$ is a non-stationary time series.
\end{itemize}
In the following, we provide an example to investigate the process.
\begin{example}
Let $(S_{t})_{t\geq0}$ be a SLCP process with period $\tau=6.5$. Furthermore, the lengths of the 
successive subintervals of each period interval are 0.5, 2.5, 3, 0.5 where corresponding occurrence rates 
of the semi-L\'evy Poisson process on these subintervals are assumed as 4, 10, 5, 30. 
Moreover, the distribution of jump sizes on these subintervals are N(2, 4), N(1.5, 2.5), N(2.5, 1.5) and N(1.75, 3), 
where $N(\mu, \sigma^{2})$ denotes a Normal distribution with mean $\mu$ and variance $\sigma^{2}$.
In this example, we consider SLD-COGARCH(1,3) process with parameters $\alpha_{0}=10^{-6}$, $\alpha_{1}=0.005$, $\beta_{1}=2.1$, 
$\beta_{2}=6$ and $\beta_{3}=0.6$. Thus, the matrix $B$ is
\begin{align*}
B = \begin{pmatrix}
0 & 1 & 0 \\
0 & 0 & 1 \\
-0.6 & -6 & -2.1\\
\end{pmatrix}
\end{align*} 
and conditions (\ref{equ3.2}), (\ref{equ3.4}) and (\ref{equ3.5}) are satisfied. For such SLD-COGARCH process, we simulate 
$G_{\Upsilon_{n}}$ 
for the duration of $m=30$ period intervals by the specified parameters, $\mathbf{Y}_{0}=( 0.37\times10^{-3}, 0.05\times10^{-3}, 
0.19\times10^{-3})^{\prime}$ and $G_{0}=0$ by the suggested simulation algorithm. 
Then, by the step 6 of the last algorithm, we get equally space samples with length $l=0.25$. So, we have 780 discretized 
samples of this 30 period intervals.
\end{example}
\begin{figure}[H]
\centering
\includegraphics[scale=0.35]{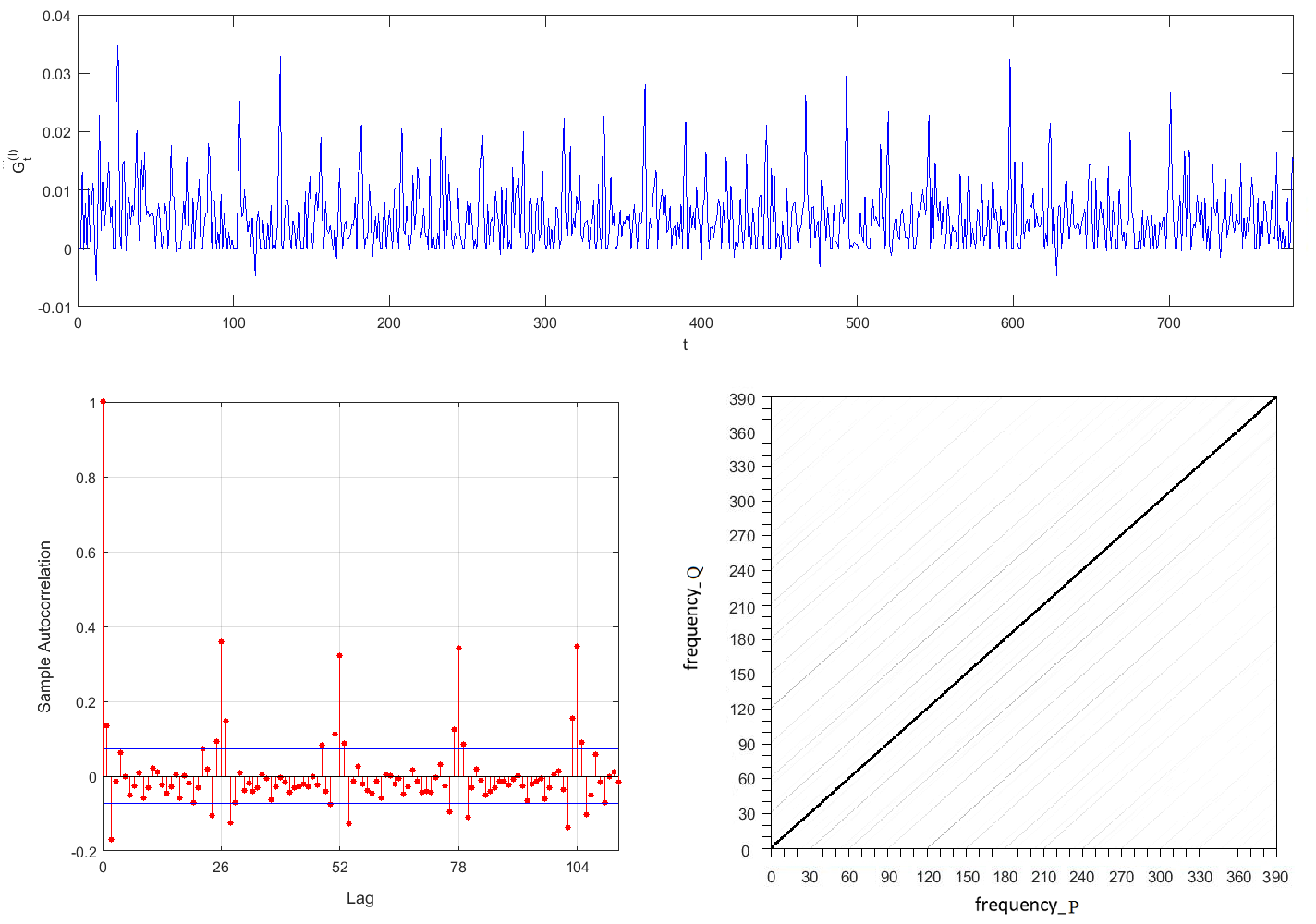}
\caption{\scriptsize{Top: the increments ($G_{il}^{(l)}=G_{(i+1)l}-G_{il}$) of the sampled 
process $(G_{il})_{i=0,\cdots,779}$; bottom left: the sample autocorrelation plot of $(G_{il}^{(l)})_{i=0,\cdots,779}$; 
bottom right: the significant values of the sample spectral coherence statistic of 
$(G_{il}^{(l)})_{i=0,\cdots,779}$ with $\alpha=0.05$ and $M=240$.}}
\end{figure}
In Figure 1, we see the differenced series $(G_{il}^{(l)})_{i=0, 1, \cdots, 779}$ (top) and its sample autocorrelation (bottom left). The sample spectral coherence statistic of this series for a specified collection of pairs $(P,Q)$ and $M=240$ that exceed the threshold corresponding to $\alpha=0.05$ which is presented at the bottom right. If $X_{i}:=G_{il}^{(l)}$, then the parallel lines for the sample spectral coherence confirm the series $(X_{i})_{i=0, 1, \cdots, 779}$ is PC with period $\varrho=780/30=26$. This implies that, for $i, j=0, 1, \cdots, 779$, 
\begin{align}\label{equ4.6}
E\big(G_{il}^{(l)}\big)=E(X_{i})=E(X_{i+\varrho})=E\big(G_{il+\varrho l}^{(l)}\big)
\end{align}
and
\begin{align}\label{equ4.7}
cov\big(G_{il}^{(l)},G_{jl}^{(l)}\big)=cov(X_{i},X_{j})=cov(X_{i+\varrho},X_{j+\varrho})=cov\big(G_{il+\varrho l}^{(l)},G_{jl+\varrho 
l}^{(l)}\big).
\end{align}

So, it follows from (\ref{equ4.6}) and (\ref{equ4.7}) that the series $(G_{t}^{(l)})_{t=0, l, \cdots, 779l}$ has a PC structure with period $\tau=\varrho l=6.5$. 

\section{Real data analysis}
\setcounter{equation}{0}
One of the striking features of financial time series is that the log returns have negligible correlation while its square log returns are significantly correlated \cite{b5}. Many of these time series are collected at a high-frequency and typically show a PC structure in their squared log returns \cite{r2}. In this section, we evaluate the results of Corollary 3.4 and show that the SLD-COGARCH process can capture the periodic structure of high-frequency financial data. For this, we use the 15-minute log returns of the Dow Jones Industrial Average (DJIA) index.

This data set is recorded between 9:35 to 16:00 from September 4th, 2015 to August 14th, 2018. There was a total of 735 trading days not including the weekends and holidays with 26 15-minute observations per day, resulting in the total of n = 19110 15-minute observations. Figure 2, shows the DJIA index for the specified times.
\begin{figure}[H]
\centering
\includegraphics[width=10.5cm,height=7.5cm]{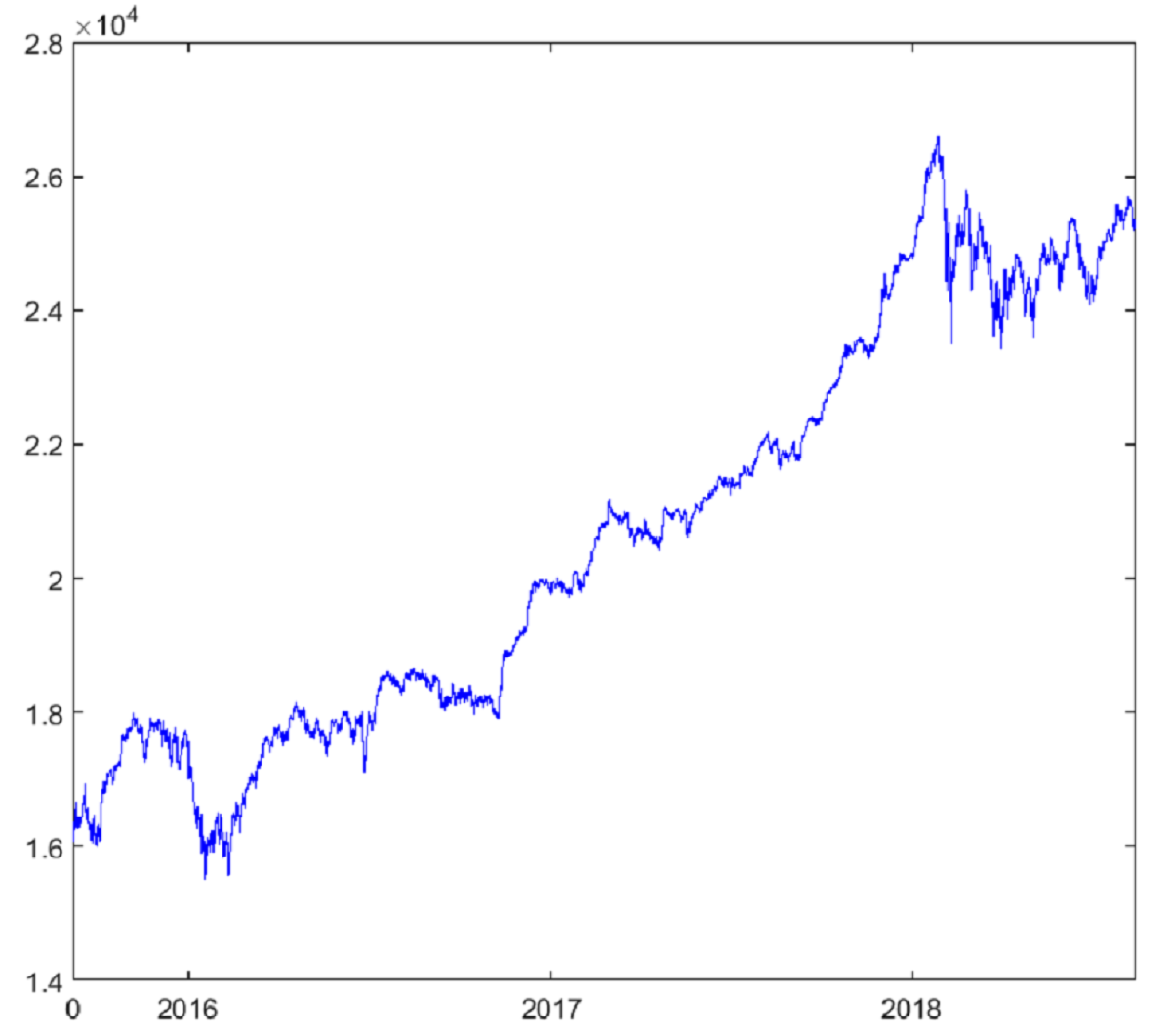}
\caption{\scriptsize{DJIA 15-minute data from September 4th, 2015 to August 14th, 2018.}}
\end{figure}
\begin{figure}[H]
\centering
\includegraphics[scale=0.34]{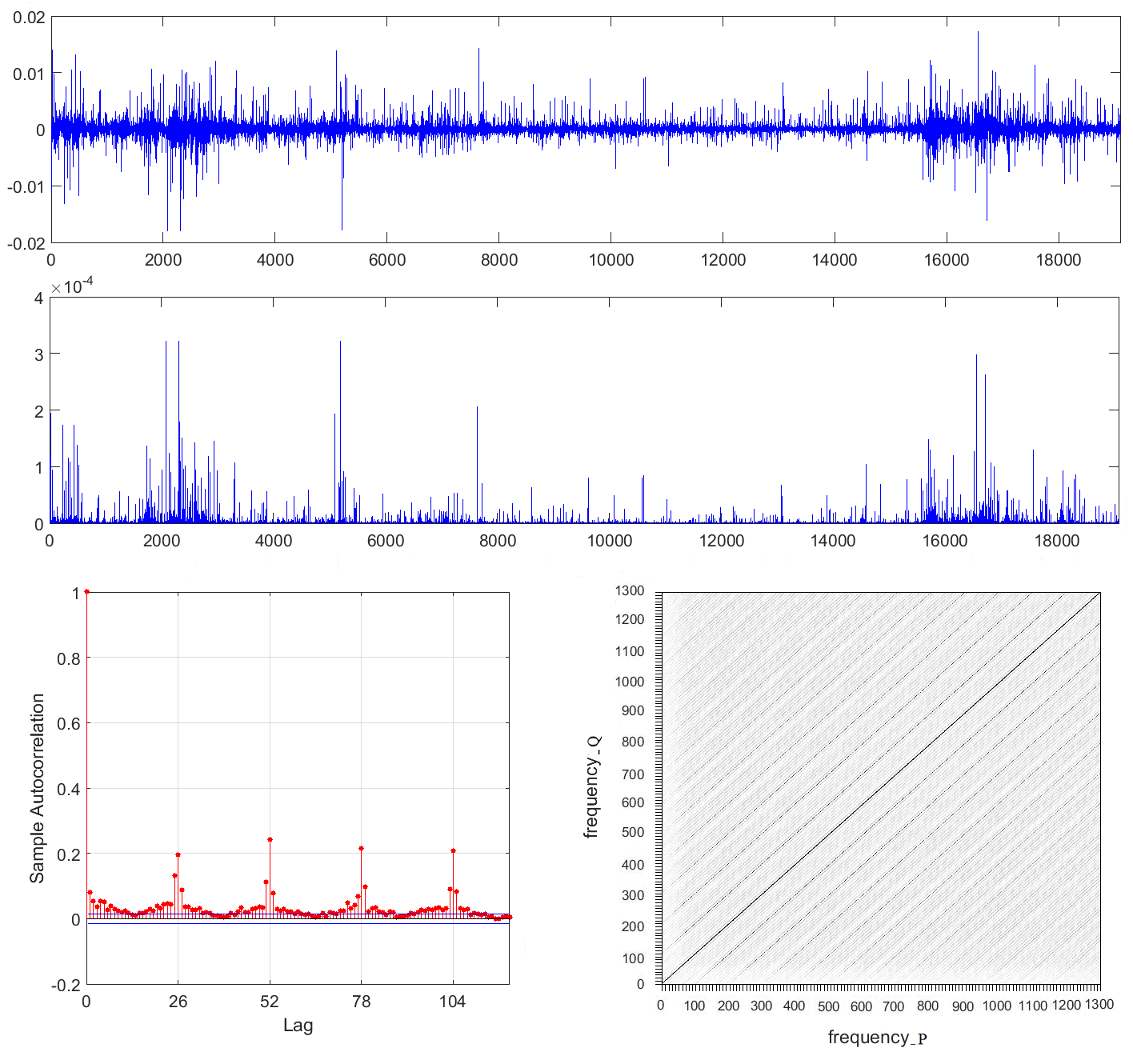}
\caption{\scriptsize{Top: DJIA log returns; middle: their squares; bottom left: the sample
autocorrelation plot of the squared log returns and bottom right: the significant values of the sample spectral coherence statistic of the squared log 
returns.}}
\end{figure} 
In Figure 3, we show the 15-minute log returns of DJIA (top), their squares (middle) and sample autocorrelation function (ACF) of the squared log returns (bottom left). It is clear from ACF that 15-minute squared log returns have a seasonal structure with period 26. To show more precisely the PC structure of these data, we sample last 100 trading days. These data were recorded from March 22th, 2018 to August 14th, 2018 which contain 2600 observations. Then, we present the sample spectral coherence statistic of their squared log returns for a specified collection of pairs $(P,Q)$ and $M=100$. It is shown that it exceeds the threshold corresponding to $\alpha=0.05$ in Figure 3 (bottom right). Therefore, sample spectral coherence statistic clearly shows that the squared log returns have a PC structure with period $\varrho=2600/100=26$.

Currently we show that how well SLD-COGARCH process fits the 15-minute log returns of DJIA index. For this, we first simulate the driving process $(S_{t})_{t\geq0}$ which is specified in Example 4.1 in the case that $\delta=0$ and the distribution of jump sizes are assumed to be N(0,4), N(0,2.5), N(0,1.5) and N(0,3). Second, we simulate the SLD-COGARCH(1,3) process for the duration of $m=735$ period intervals with the initial value $\mathbf{Y}_{0}=( 0.37\times10^{-3}, 0.05\times10^{-3}, 0.19\times10^{-3})^{\prime}$ and the parameters $\alpha_{0}=0.8\times10^{-6}$, $\alpha_{1}=0.0275$, $\beta_{1}=2.1$, $\beta_{2}=6$ and $\beta_{3}=0.6$ by the suggested simulation algorithm. Then, by the step 6 of the last algorithm, we get 19110 equally space samples by length $l=0.25$. 
\begin{figure}[H]
\centering
\includegraphics[scale=0.34]{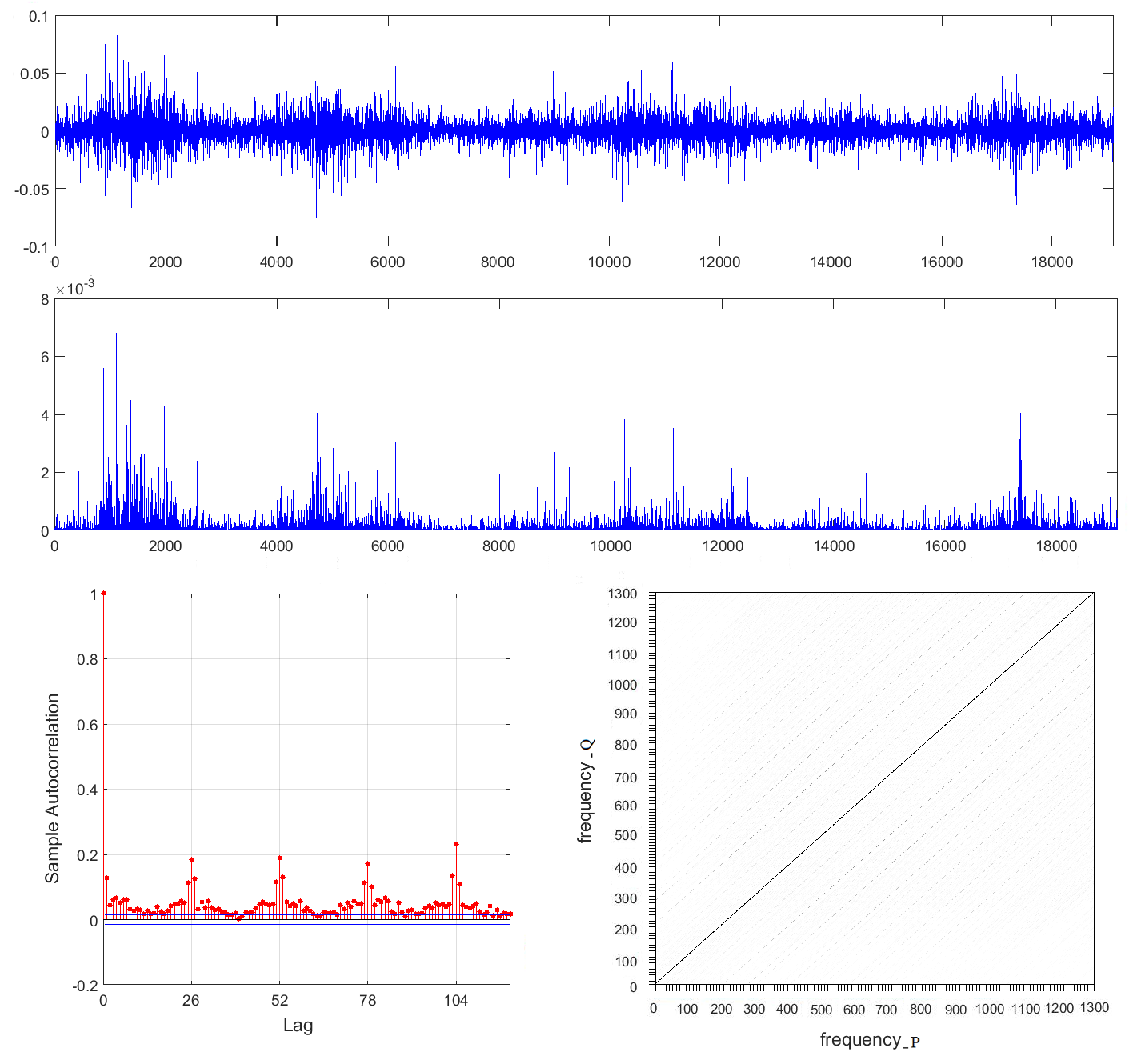}
\caption{\scriptsize{Top: the increments $(G_{il}^{(l)}=G_{(i+1)l}-G_{il})$ of the sampled process $(G_{il})_{i=0,\cdots,19109}$; middle: their squares; bottom left: the sample autocorrelation plot of the squared increments; bottom right: the significant values of the sample spectral coherence statistic of the squared increments with $\alpha=0.05$ and $M=550$.}}
\end{figure}
Figure 4 shows the increments of the sampled process $(G_{il})_{i=0,\cdots,19109}$, (top) and their squares (middle) and ACF of the squared increments (bottom left). We select the last 2600 simulated data as sample and present the sample spectral coherence statistic of their squared increments for a specified collection of pairs $(P,Q)$ and $M=550$. It is shown that it exceeds the threshold corresponding to $\alpha=0.05$ in Figure 4 (bottom right). The ACF and sample spectral coherence statistic demonstrate that the squared increments have a PC structure with period $\varrho=2600/100=26$. Therefore, these graphs show that SLD-COGARCH process can be consistent with 15-minute log returns of DJIA index.

\section{Appendix}
\setcounter{equation}{0}
\textbf{P1: Proof of Lemma 2.1}\\
For any $t\geq0$ there exist $r=1, \cdots, d$ such that $t\in A_{r+md}=(s_{_{r-1+md}}, s_{_{r+md}}]$, where $m=[\frac{t}{\tau}]$. 
So, it follows from Definition 2.1 and Definition 2.2 that
\begin{align*}
E\big(e^{iuS_{t}}\big)&=E\big(e^{iu\big(\delta t+\sum_{n=1}^{N(t)}Z_{n}\big)}\big)
=e^{iu\delta t}E\big(e^{iu\sum_{n=1}^{N(t)}Z_{n}}\big)\\
&=e^{iu\delta t}E\Big(e^{iu\sum_{n=1}^{N(s_{_{1}})}Z_{n}^{(1)}}\times 
e^{iu\sum_{n=N(s_{_{1}})+1}^{N(s_{_{2}})}Z_{n}^{(2)}}\times\cdots\times 
e^{iu\sum_{n=N(s_{_{r-1+md}})+1}^{N(t)}Z_{n}^{(r)}}\Big)\\
&=e^{iu\delta t}E\Big(\prod_{j=1}^{d}\prod_{k=0}^{m-1}
e^{iu\sum_{n=N(s_{_{j-1+kd}})+1}^{N(s_{_{j+kd}})}Z_{n}^{(j)}}
\prod_{j=1}^{r-1}e^{iu\sum_{n=N(s_{_{j-1+md}})+1}^{N(s_{_{j+md}})}Z_{n}^{(j)}}\\
&\qquad\qquad\times e^{iu\sum_{n=N(s_{_{r-1+md}})+1}^{N(t)}Z_{n}^{(r)}}\Big).
\end{align*}
Since, by Definition 2.2, the jump sizes $Z_{n}^{(j)}, n\in\mathbb{N},$ are independent for $j=1,\cdots,d$, we have that
\begin{align}\nonumber
E\big(e^{iuS_{t}}\big)&=e^{iu\delta t}\prod_{j=1}^{d}\prod_{k=0}^{m-1}
E\big(e^{iu\sum_{n=N(s_{_{j-1+kd}})+1}^{N(s_{_{j+kd}})}Z_{n}^{(j)}}\big)
\prod_{j=1}^{r-1}E\big(e^{iu\sum_{n=N(s_{_{j-1+md}})+1}^{N(s_{_{j+md}})}Z_{n}^{(j)}}\big)\\ \label{Lem2.1(5.1)}
&\quad\times E\big(e^{iu\sum_{n=N(s_{_{r-1+md}})+1}^{N(t)}Z_{n}^{(r)}}\big).
\end{align}
By the fact that the jump sizes $Z_{n}^{(j)}, n\in\mathbb{N},$ in each partition $A_{j+kd}=(s_{_{j-1+kd}},s_{_{j+kd}}], 
k\in\mathbb{N}^{0},$ 
are i.i.d. random variables with distribution $F_{j}$, it follows from Definition 2.1 and conditional expected value that
\begin{align}\nonumber
&E\big(e^{iu\sum_{n=N(s_{_{j-1+kd}})+1}^{N(s_{_{j+kd}})}Z_{n}^{(j)}}\big)=E\Big(E\big(e^{iu\sum_{n=1}^{N(s_{_{j+kd}})-
N(s_{_{j-1+kd}})}Z_{n}^{(j)}}\vert N(s_{_{j+kd}})-N(s_{_{j-1+kd}})=N\big)\Big)\\ \nonumber
&\qquad\qquad=\sum_{N=0}^{\infty}E\big(e^{iu\sum_{n=1}^{N}Z_{n}^{(j)}}\big)P\big(N(s_{_{j+kd}})-N(s_{_{j-1+kd}})=N\big)\\ 
\nonumber
&\qquad\qquad=\sum_{N=0}^{\infty}\big(E(e^{iuZ_{n}^{(j)}})\big)^{N}\frac{\big(\Lambda(s_{_{j+kd}})-\Lambda(s_{_{j-1+kd}})\big)^{N}
e^{-\big(\Lambda(s_{_{j+kd}})-\Lambda(s_{_{j-1+kd}})\big)}}{N!}\\ \nonumber
&\qquad\qquad=e^{-(\Lambda(s_{_{j+kd}})-\Lambda(s_{_{j-1+kd}}))}\sum_{N=0}^{\infty}\frac{\Big(\int_{\mathbb{R}}e^{iuz}
\big(\Lambda(s_{_{j+kd}})-\Lambda(s_{_{j-1+kd}})\big)F_{j}(dz)
\Big)^{N}}{N!}\\ \nonumber
&\qquad\qquad=e^{-(\Lambda(s_{_{j+kd}})-\Lambda(s_{_{j-1+kd}}))}e^{\int_{\mathbb{R}}e^{iuz}
(\Lambda(s_{_{j+kd}})-\Lambda(s_{_{j-1+kd}}))F_{j}(dz)}\\ \label{Lem2.1(5.2)}
&\qquad\qquad=e^{\int_{\mathbb{R}}(e^{iuz}-1)\lambda_{j}|A_{j}|F_{j}(dz)},
\end{align}
where $|A_{j}|=|A_{j+kd}|=s_{_{j}}-s_{_{j-1}}$. The last equality follows from (\ref{equ2.1}) and the fact that the semi-L\`evy Poisson 
process 
$N(t)$ has the occurrence rate $\lambda_{j}$ 
on $A_{j+kd}$ for $k\in\mathbb{N}^{0}$. An argument similar to that used in (\ref{Lem2.1(5.2)}) shows that 
\begin{align}\label{Lem2.1(5.3)}
E\big(e^{iu\sum_{n=N(s_{_{r-1+md}})+1}^{N(t)}Z_{n}^{(r)}}\big)=e^{\int_{\mathbb{R}}(e^{iuz}-1)\lambda_{r}(t-s_{_{r-1+md}})F_{r}
(dz)}
\end{align}
By replacing from (\ref{Lem2.1(5.2)}) and (\ref{Lem2.1(5.3)}) in (\ref{Lem2.1(5.1)}), it reads as
\begin{align*}
E\big(e^{iuS_{t}}\big)&=e^{iu\delta t}\prod_{j=1}^{d}e^{\int_{\mathbb{R}}(e^{iuz}-1)m\lambda_{j}|A_{j}|F_{j}(dz)}
\prod_{j=1}^{r-1}e^{\int_{\mathbb{R}}(e^{iuz}-1)\lambda_{j}|A_{j}|F_{j}(dz)}\\
&\quad\times e^{\int_{\mathbb{R}}(e^{iuz}-1)\lambda_{r}(t-s_{_{r-1+md}})F_{r}(dz)}\\
&=e^{\Big(iu\big(\delta t+\int_{|z|\leq1}z\nu_{t}(dz)\big)+\int_{\mathbb{R}}(e^{iuz}-1-iuzI_{|z|\leq1})\nu_{t}(dz)\Big)},
\end{align*}
where
\begin{align*}
\nu_{t}(dz)=\sum_{j=1}^{r-1}(m+1)\lambda_{j}|A_{j}|F_{j}(dz)+\sum_{j=r}^{d}m\lambda_{j}|A_{j}|F_{j}(dz)+
\lambda_{r}(t-s_{_{r-1+md}})F_{r}(dz).
\end{align*}
\\
\textbf{P2: Proof of Theorem 3.1}
\vspace{2mm}

\textbf{Proof ($a$):} Let $(S_{t})_{t\geq0}$ be the driving process introduced by Definition 2.2, $\Upsilon_{_{n}}$ be the time of $n^{th}$ jump and $\Upsilon_{_{0}}=0$. Following the method of Brockwell et al. \cite{b5}, and as $(S_{t})_{t\geq0}$  has no jump over $[\Upsilon_{_{n}},\Upsilon_{_{n+1}})$, it follows by (\ref{equ2.5}) that $[S,S]_{t}=0$ for $t\in[\Upsilon_{_{n}},\Upsilon_{_{n+1}})$. So (\ref{equ2.7}) implies that $\mathbf{Y}_{t}$ satisfies $d\mathbf{Y}_{t}=B\mathbf{Y}_{t}dt$ for $t\in[\Upsilon_{_{n}},\Upsilon_{_{n+1}})$.  Hence
\begin{align*}
e^{-Bt}d\mathbf{Y}_{t}=Be^{-Bt}\mathbf{Y}_{t}dt
\end{align*}
and
\begin{align}\label{Th3.1(5.1)}
d\big(e^{-Bt}\mathbf{Y}_{t}\big)=0.
\end{align}
Integrating from $\Upsilon_{_{n}}$ to $t\in[\Upsilon_{_{n}},\Upsilon_{_{n+1}})$, we have that
\begin{align}\label{Th3.1(5.2)}
\mathbf{Y}_{t}=e^{B(t-\Upsilon_{_{n}})}\mathbf{Y}_{\Upsilon_{_{n}}}.
\end{align}
By (\ref{equ2.7}) and (\ref{equ2.5}), we have that $\mathbf{Y}_{t}$ has a 
jump of size $\mathbf{e}(\alpha_{0}+\mathbf{a}^{\prime}\mathbf{Y}_{\Upsilon_{_{n}}^{-}})Z_{_{n}}^{2}$ 
at $\Upsilon_{_{n}}$. So, $\mathbf{Y}_{t}$ at $\Upsilon_{_{n}}$ can be written as
\begin{align}\nonumber
\mathbf{Y}_{\Upsilon_{_{n}}}&=\mathbf{Y}_{\Upsilon_{_{n}}^{-}}+\mathbf{e}
(\alpha_{0}+\mathbf{a}^{\prime}\mathbf{Y}_{\Upsilon_{_{n}}^{-}})Z_{_{n}}^{2}\\ \label{Th3.1(5.3)}
&=(I+Z_{_{n}}^{2}\mathbf{e}\mathbf{a}^{\prime})\mathbf{Y}_{\Upsilon_{_{n}}^{-}}+
\alpha_{0}Z_{_{n}}^{2}\mathbf{e}.
\end{align}
Integrating of (\ref{Th3.1(5.1)}) this time from $\Upsilon_{_{n-1}}$ to $\Upsilon_{_{n}}^{-}$, we have that 
\begin{align}\label{Th3.1(5.4)}
\mathbf{Y}_{\Upsilon_{_{n}}^{-}}=e^{B(\Upsilon_{_{n}}-\Upsilon_{_{n-1}})}\mathbf{Y}_{\Upsilon_{_{n-1}}}.
\end{align}
Therefore 
\begin{align}\label{Th3.1(5.5)}
\mathbf{Y}_{\Upsilon_{_{n}}}=(I+Z_{_{n}}^{2}\mathbf{e}\mathbf{a}^{\prime})
e^{B(\Upsilon_{_{n}}-\Upsilon_{_{n-1}})}\mathbf{Y}_{\Upsilon_{_{n-1}}}+\alpha_{0}Z_{_{n}}^{2}\mathbf{e}.
\end{align}
By assuming $t\in[\Upsilon_{_{N(t)}},\Upsilon_{_{N(t)+1}})$, (\ref{Th3.1(5.2)}) reads as
\begin{align}\label{Th3.1(5.6)}
\mathbf{Y}_{t}=e^{B(t-\Upsilon_{_{N(t)}})}\mathbf{Y}_{\Upsilon_{_{N(t)}}}
\end{align}
and (\ref{Th3.1(5.5)}) as
\begin{align}\label{Th3.1(5.7)}
\mathbf{Y}_{\Upsilon_{_{N(t)}}}=(I+Z_{_{N(t)}}^{2}\mathbf{e}\mathbf{a}^{\prime})
e^{B(\Upsilon_{_{N(t)}}-\Upsilon_{_{N(t)-1}})}\mathbf{Y}_{\Upsilon_{_{N(t)-1}}}
+\alpha_{0}Z_{_{N(t)}}^{2}\mathbf{e}.
\end{align}
By successive replacement from (\ref{Th3.1(5.7)}) in (\ref{Th3.1(5.6)}) we have that
\begin{align}\nonumber
\mathbf{Y}_{t}&=e^{B\big(t-\Upsilon_{_{N(t)}}\big)}\Big[\big(I+Z_{_{N(t)}}^{2}\mathbf{e}\mathbf{a}^{\prime}\big)
e^{B\big(\Upsilon_{_{N(t)}}-\Upsilon_{_{N(t)-1}}\big)}\times\\ \nonumber
&\qquad\qquad\qquad\cdots\times\big(I+Z_{_{N(s)+2}}^{2}\mathbf{e}\mathbf{a}^{\prime}\big)
e^{B\big(\Upsilon_{_{N(s)+2}}-\Upsilon_{_{N(s)+1}}\big)}\Big]\mathbf{Y}_{\Upsilon_{_{N(s)+1}}}\\ \nonumber
&+e^{B\big(t-\Upsilon_{_{N(t)}}\big)}\Big[(\alpha_{0}Z_{_{N(t)}}^{2}\mathbf{e})+\sum_{i=0}^{N(t)-N(s)-3}
\big(I+Z_{_{N(t)}}^{2}\mathbf{e}\mathbf{a}^{\prime}\big)
e^{B\big(\Upsilon_{_{N(t)}}-\Upsilon_{_{N(t)-1}}\big)}\times\\ \label{Th3.1(5.8)}
&\qquad\qquad\qquad\cdots\times\big(I+Z_{_{N(t)-i}}^{2}\mathbf{e}\mathbf{a}^{\prime}\big)
e^{B\big(\Upsilon_{_{N(t)-i}}-\Upsilon_{_{N(t)-i-1}}\big)}(\alpha_{0}Z_{_{N(t)-i-1}}^{2}\mathbf{e})\Big].
\end{align}
Similar to (\ref{Th3.1(5.3)}) $\mathbf{Y}_{t}$ at $\Upsilon_{_{N(s)+1}}$ can be written as
\begin{align}\label{Th3.1(5.9)}
\mathbf{Y}_{\Upsilon_{_{N(s)+1}}}=(I+Z_{_{N(s)+1}}^{2}\mathbf{e}\mathbf{a}^{\prime})
\mathbf{Y}_{\Upsilon_{_{N(s)+1}}^{^{-}}}+\alpha_{0}Z_{_{N(s)+1}}^{2}\mathbf{e}.
\end{align}
By integrating of (\ref{Th3.1(5.1)}) from $s$ to $\Upsilon_{_{N(s)+1}}$, we have that
\begin{align*}
\mathbf{Y}_{\Upsilon_{_{N(s)+1}}^{^{-}}}=e^{B(\Upsilon_{_{N(s)+1}}-s)}\mathbf{Y}_{s}.
\end{align*}
Replacing this in (\ref{Th3.1(5.9)}) implies that
\begin{align}\label{Th3.1(5.10)}
\mathbf{Y}_{\Upsilon_{_{N(s)+1}}}=(I+Z_{_{N(s)+1}}^{2}\mathbf{e}\mathbf{a}^{\prime})
e^{B(\Upsilon_{_{N(s)+1}}-s)}\mathbf{Y}_{s}+\alpha_{0}Z_{_{N(s)+1}}^{2}\mathbf{e}.
\end{align}
By replacing from (\ref{Th3.1(5.10)}) in (\ref{Th3.1(5.8)}), it reads as
\begin{align}\nonumber
\mathbf{Y}_{t}&=e^{B\big(t-\Upsilon_{_{N(t)}}\big)}\Big[\big(I+Z_{_{N(t)}}^{2}\mathbf{e}\mathbf{a}^{\prime}\big)
e^{B\big(\Upsilon_{_{N(t)}}-\Upsilon_{_{N(t)-1}}\big)}\times\\ \nonumber
&\quad\qquad\qquad\cdots\times
\big(I+Z_{_{N(s)+2}}^{2}\mathbf{e}\mathbf{a}^{\prime}\big)
e^{B\big(\Upsilon_{_{N(s)+2}}-\Upsilon_{_{N(s)+1}}\big)}\big(I+Z_{_{N(s)+1}}^{2}\mathbf{e}\mathbf{a}^{\prime}\big)
e^{B\big(\Upsilon_{_{N(s)+1}}-s\big)}\Big]\mathbf{Y}_{s}\\ \nonumber
&+e^{B\big(t-\Upsilon_{_{N(t)}}\big)}\Big[(\alpha_{0}Z_{_{N(t)}}^{2}\mathbf{e})+\sum_{i=0}^{N(t)-N(s)-2}
\big(I+Z_{_{N(t)}}^{2}\mathbf{e}\mathbf{a}^{\prime}\big)
e^{B\big(\Upsilon_{_{N(t)}}-\Upsilon_{_{N(t)-1}}\big)}\times\\ \nonumber
&\quad\qquad\qquad\cdots\times\big(I+Z_{_{N(t)-i}}^{2}\mathbf{e}\mathbf{a}^{\prime}\big)
e^{B\big(\Upsilon_{_{N(t)-i}}-\Upsilon_{_{N(t)-i-1}}\big)}(\alpha_{0}Z_{_{N(t)-i-1}}^{2}\mathbf{e})\Big]\\ 
\label{Th3.1(5.11)}
&=:J_{s,t}\mathbf{Y}_{s}+\mathbf{K}_{s,t}.
\end{align}

\textbf{Proof (b):} We show that the distribution of $J_{s,t}$ and $\mathbf{K}_{s,t}$ are periodic in both indices $s$ and $t$, with the same period $\tau$. Following the fact that $J_{s,t}$ and $\mathbf{K}_{s,t}$ which are defined by (\ref{Th3.1(5.11)}) depend on $s$ and $t$ through
\begin{itemize}
\item
the jump sizes $Z_{\cdot}$'s which by Definition 2.1 are independent and periodically identically distributed random variables with distribution $F_{j}, j=1, \cdots, d$,
\item
$\Upsilon_{_{N(t)}}-\Upsilon_{_{N(t)-1}}, \cdots, \Upsilon_{_{N(s)+2}}-\Upsilon_{_{N(s)+1}}$ which are inter-arrival times of $N(t)$ which by Definition 2.1 has periodically stationary increment with period $\tau$, So are periodic with period $\tau$. Also $N(t)-N(s)$ are periodic by the same reason,
\item
$t-\Upsilon_{_{N(t)}}$ is the current life at time $t$ which is periodic with period $\tau$. Also $J_{s,t}$ is a function of $\Upsilon_{_{N(s)+1}}-s$ which excess life and is periodic with period $\tau$, 
\end{itemize}
it follows that 
\begin{align*}
J_{s,t}\overset{d}{=}J_{s+\tau,t+\tau} \qquad \textnormal{and}\qquad \mathbf{K}_{s,t}\overset{d}{=}\mathbf{K}_{s+\tau,t+\tau}.
\end{align*}
So, $(J_{s,t},\mathbf{K}_{s,t})$ are periodic in both indices $s$ and $t$, with the same period $\tau$.

\vspace{5mm}
\textbf{Proof (c):} By (\ref{Th3.1(5.11)}) it is clear that 
$(J_{s,t},\mathbf{K}_{s,t})$ and $(J_{r,u},\mathbf{K}_{r,u})$ are constructed from segments of the driving process $(S_{t})_{t\geq0}$ in the intervals $(s,t]$ and $(r,u]$, respectively. Thus $(J_{s,t},\mathbf{K}_{s,t})$ and $(J_{r,u},\mathbf{K}_{r,u})$ are independent for $0\leq s\leq t\leq r\leq u$.\\
Now we show that $\mathbf{Y}_{s}$ is independent of $\big(J_{s,t},\mathbf{K}_{s,t}\big)$. By (\ref{Th3.1(5.11)}) we have that
\begin{align*}
\mathbf{Y}_{s}=J_{0,s}\mathbf{Y}_{0}+\mathbf{K}_{0,s}.
\end{align*}
As $(J_{0,s},\mathbf{K}_{0,s})$ and $(J_{s,t},\mathbf{K}_{s,t})$ are independent and that $\mathbf{Y}_{0}$ is independent of the driving process $(S_{t})_{t\geq0}$, assumed by (\ref{equ2.7}), it follows that $\mathbf{Y}_{s}$ and $(J_{s,t},\mathbf{K}_{s,t})$ are independent.
\vspace{6mm}
\\
\textbf{P3: Proof of Theorem 3.2}\\
For fixed $t\in[0,\tau)$ and for $m\in\mathbb{N}$, it follows from (\ref{Th3.1(5.11)}) that 
\begin{align}\label{Th3.2(5.12)}
\mathbf{Y}_{t+m\tau}=J_{t+(m-1)\tau,t+m\tau}\mathbf{Y}_{t+(m-1)\tau}+\mathbf{K}_{t+(m-1)\tau,t+m\tau}.
\end{align}
By successive replacement from (\ref{Th3.2(5.12)}) we obtain
\begin{align}\nonumber
\mathbf{Y}_{t+m\tau}&=\Big(J_{t+(m-1)\tau,t+m\tau}\cdots J_{t,t+\tau}\Big)\mathbf{Y}_{t}
+\Big[\mathbf{K}_{t+(m-1)\tau,t+m\tau}\\ \label{Th3.2(5.13)}
&\quad +J_{t+(m-1)\tau,t+m\tau}\mathbf{K}_{t+(m-2)\tau,t+(m-1)\tau}+\cdots+ J_{t+(m-1)\tau,t+m\tau}\cdots J_{t+\tau,t+2\tau}\mathbf{K}_{t,t+\tau}\Big].
\end{align}
Since, by the part ($b$) of Theorem 3.1, $(J_{t,t+\tau}, \mathbf{K}_{t,t+\tau}), (J_{t+\tau,t+2\tau}, \mathbf{K}_{t+\tau,t+2\tau}),\cdots$ are i.i.d., it follows from (\ref{Th3.2(5.13)}) that
\begin{align}\nonumber
\mathbf{Y}_{t+m\tau}&\overset{d}{=}\Big(J_{t,t+\tau}\cdots J_{t+(m-1)\tau,t+m\tau}\Big)\mathbf{Y}_{t}
+\Big[\mathbf{K}_{t,t+\tau}+J_{t,t+\tau}\mathbf{K}_{t+\tau,t+2\tau}\\ \nonumber
&\hspace{5cm}+\cdots+ J_{t,t+\tau}\cdots J_{t+(m-2)\tau,t+(m-1)\tau}\mathbf{K}_{t+(m-1)\tau,t+m\tau}\Big]\\ \nonumber
&=\Big(\prod_{i=1}^{m}J_{t+(i-1)\tau, t+i\tau}\Big)\mathbf{Y}_{t}
+\Big[\mathbf{K}_{t,t+\tau}+\sum_{i=1}^{m-1}J_{t,t+\tau}\cdots J_{t+(i-1)\tau,t+i\tau}
\mathbf{K}_{t+i\tau,t+(i+1)\tau}\Big]\\ \label{Th3.2(5.14)}
&=:H_{m}^{(t)}\mathbf{Y}_{t}+\mathbf{U}_{m}^{(t)}.
\end{align}
The following method of Brockwell et al. \cite{b5}, for the proof of this theorem we use the general theory of multivariate random recurrence equations, see Bougerol and Picard \cite{b3}. For this, we first show that there exist 
some $r\in[1,\infty]$ and $(q\times q)-$matrix $P$ such that
\begin{align*}
E\Big(log||J_{t,t+\tau}||_{P,r}\Big)<0 \qquad \textnormal{and}\qquad E\Big(log^{+}||\mathbf{K}_{t,t+\tau}||
_{P,r}\Big)<\infty,
\end{align*}
where $log^{+}(x)=log(max\lbrace1,x\rbrace)$. Then, it follows from the strong law of large numbers that the Lyapunov exponent of the i.i.d. sequence $(J_{t+(n-1)\tau,t+n\tau})_{n\in\mathbb{N}}$ is strictly negative
\begin{align*}
\gamma&=\lim_{n\rightarrow\infty}\frac{1}{n}log||J_{t,t+\tau}\cdots J_{t+(n-1)\tau,t+n\tau}||_{P,r}\\
&\leq \lim_{n\rightarrow\infty}\frac{1}{n}\sum_{i=1}^{n}log||J_{t+(i-1)\tau,t+i\tau}||_{P,r}<0
\end{align*}
Therefore, as shown by Bougerol and Picard \cite{b3}, the sum $\mathbf{U}_{m}^{(t)}$ for all $t\in[0,\tau)$ is converges almost surely (a.s.).

Since the eigenvalues of the matrix $B$, $\eta_{1}, \cdots, \eta_{q}$, are distinct, we define the matrix $P$ as follows
\begin{align}\label{Th3.2(5.15)}
P:= \begin{pmatrix}
1 & \cdots & 1\\
\eta_{1} & \cdots & \eta_{q}\\
\vdots & \cdots & \vdots \\
\eta_{1}^{q-1} & \cdots & \eta_{q}^{q-1}
\end{pmatrix}.
\end{align}
So, the matrix $C:=P^{-1}BP$ is the diagonal matrix with entries $\eta_{1}, \cdots, \eta_{q}$ on the diagonal. Using the fact that $P$ is a invertible matrix, it follows from the exponential matrix definition that, for $t\geq0$, 
\begin{align*}
e^{Ct}=e^{P^{-1}BtP}&=I+P^{-1}BtP+\frac{(P^{-1}BtP)^{2}}{2!}+\cdots\\
&=P^{-1}\big[I+Bt+\frac{(Bt)^{2}}{2!}+\cdots\big]P=P^{-1}e^{Bt}P.
\end{align*}
Hence, $P^{-1}e^{Bt}P=e^{Ct}$ is the diagonal matrix with entries $e^{\eta_{1}t}, \cdots, e^{\eta_{q}t}$ on the diagonal. Since the eigenvalues of the matrix B have strictly negative real parts, it follows from the Definition 2.5 that, for $r\in[1,\infty]$, 
\begin{align}\label{Th3.2(5.16)}
||e^{Bt}||_{P,r}=||Pe^{Ct}P^{-1}||_{P,r}=||P^{-1}Pe^{Ct}P^{-1}P||_{r}=||e^{Ct}||_{r}\leq e^{\eta t}.
\end{align}
By denoting
\begin{align}\nonumber
C_{i}&:=(I+Z_{i}^{2}\mathbf{e}\mathbf{a}^{\prime})e^{B(\Upsilon_{_{i}}-\Upsilon_{_{i-1}})},\\ \label{Th3.2(5.17)}
\mathbf{D}_{i}&:=\alpha_{0}Z_{i}^{2}\mathbf{e},\qquad i\in\mathbb{N},
\end{align}
in (\ref{Th3.1(5.11)}) we have that
\begin{align*}
J_{t,t+\tau}&=e^{B(t+\tau-\Upsilon_{_{N(t+\tau)}})}C_{_{N(t+\tau)}}\cdots C_{_{N(t)+2}}
(I+Z_{_{N(t)+1}}^{2}\mathbf{e}\mathbf{a}^{\prime})e^{B(\Upsilon_{_{N(t)+1}}-t)}\\
&=e^{B(t+\tau-\Upsilon_{_{N(t+\tau)}})}C_{_{N(t+\tau)}}\cdots C_{_{N(t)+2}}C_{_{N(t)+1}}e^{B(\Upsilon_{_{N(t)}}-t)}\\
&=e^{B(t+\tau-\Upsilon_{_{N(t+\tau)}})}\Big[\prod_{i=0}^{N(t+\tau)-N(t)-1}C_{_{N(t+\tau)-i}}\Big]e^{B(\Upsilon_{_{N(t)}}-t)}.
\end{align*}
So, by (\ref{Th3.2(5.16)}) we have that
\begin{align}\label{Th3.2(5.18)}
||J_{t,t+\tau}||_{P,r}\leq e^{\eta(t+\tau-\Upsilon_{_{N(t+\tau)}})}\Big[\prod_{i=0}^{N(t+\tau)-N(t)-1}
||C_{_{N(t+\tau)-i}}||_{P,r}\Big]e^{\eta(\Upsilon_{_{N(t)}}-t)}.
\end{align}
As (\ref{Th3.2(5.16)}) and (\ref{Th3.2(5.17)}) leads to
\begin{align}\nonumber
||C_{i}||_{P,r}&\leq (1+Z_{i}^{2}||\mathbf{e}\mathbf{a}^{\prime}||_{P,r})e^{\eta(\Upsilon_{i}-\Upsilon_{i-1})}\\\label{Th3.2(5.19)}
||\mathbf{D}_{i}||_{P,r}&=\alpha_{0}Z_{i}^{2}||\mathbf{e}||_{P,r},
\end{align}
so, (\ref{Th3.2(5.18)}) reads as
\begin{align*}
||J_{t,t+\tau}||_{P,r}&\leq e^{\eta(t+\tau-\Upsilon_{_{N(t+\tau)}})}\Big[\prod_{i=0}^{N(t+\tau)-N(t)-1}
(1+Z_{_{N(t+\tau)-i}}^{2}
||\mathbf{e}\mathbf{a}^{\prime}||_{P,r})e^{\eta(\Upsilon_{_{N(t+\tau)-i}}-\Upsilon_{_{N(t+\tau)-i-1})}}\Big]\\
&\quad\times e^{\eta(\Upsilon_{_{N(t)}}-t)}.
\end{align*}
Therefore
\begin{align}\label{Th3.2(5.20)}
log||J_{t,t+\tau}||_{P,r}\leq \eta\tau+\sum_{i=0}^{N(t+\tau)-N(t)-1}log(1+Z_{_{N(t+\tau)-i}}^{2}||\mathbf{e}\mathbf{a}^{\prime}||_{P,r}).
\end{align}
So, it follows from (\ref{equ3.2}) and (\ref{Th3.2(5.20)}) that
\begin{align*}
&E\Big(log||J_{t,t+\tau}||_{P,r}\Big)\leq \eta\tau+\Big[\sum_{n,m=0}^{\infty}\sum_{i=0}^{n-1}E\Big(log\big(1+Z_{m-i}^{2}
||\mathbf{e}\mathbf{a}^{\prime}||_{P,r}\big)\Big)\\
&\hspace{5cm}\times P\big(N(t+\tau)-N(t)=n, N(t+\tau)=m\big)\Big]\\
&\qquad<\eta\tau+\sum_{n,m=0}^{\infty}n\times\frac{-\eta\tau}{\Lambda(t+\tau)-\Lambda(t)}P\big(N(t+\tau)-N(t)=n, N(t+\tau)=m\big)\\
&\qquad=\eta\tau+\big(\Lambda(t+\tau)-\Lambda(t)\big)\frac{-\eta\tau}{\Lambda(t+\tau)-\Lambda(t)}=0
\end{align*}
To show that $E\big(log^{+}||\mathbf{K}_{t,t+\tau}||_{P,r}\big)<\infty$, observe from (\ref{Th3.1(5.11)}) and (\ref{Th3.2(5.17)}) that
\begin{align*}
\mathbf{K}_{t,t+\tau}=e^{B(t+\tau-\Upsilon_{_{N(t+\tau)}})}\Big[\mathbf{D}_{_{N(t+\tau)}}+\sum_{i=0}^{N(t+\tau)-N(t)-2}
C_{_{N(t+\tau)}}\cdots C_{_{N(t+\tau)-i}}\mathbf{D}_{_{N(t+\tau)-i-1}}\Big].
\end{align*}
So, by (\ref{Th3.2(5.16)}) we have that
\begin{align}\nonumber
||\mathbf{K}_{t,t+\tau}||_{P,r}\leq e^{\eta(t+\tau-\Upsilon_{_{N(t+\tau)}})}&\Big[||\mathbf{D}_{_{N(t+\tau)}}||_{P,r}+
\Big(\sum_{i=0}^{N(t+\tau)-N(t)-2}||C_{_{N(t+\tau)}}||_{P,r}\cdots||C_{_{N(t+\tau)-i}}||_{p,r}\\\label{Th3.2(5.21)}
&\quad\times||\mathbf{D}_{_{N(t+\tau)-i-1}}||_{P,r}\Big)\Big].
\end{align}
The eigenvalues of the matrix B have strictly negative real parts, say $\mathcal{R}e(\eta_{i})<0$. So, it follows that $\eta:=\max_{i=1,\cdots,q}\mathcal{R}e(\eta_{i})<0$ and for $t\geq0$ 
\begin{align*}
e^{\eta t}\leq1.
\end{align*}
From this and (\ref{Th3.2(5.19)}) we have that
\begin{align*}
||C_{i}||_{P,r}\leq(1+Z_{i}^{2}||\mathbf{e}\mathbf{a}^{\prime}||_{P,r}).
\end{align*}
Thus, (\ref{Th3.2(5.21)}) reads as
\begin{align*}
||\mathbf{K}_{t,t+\tau}||_{P,r}\leq \alpha_{0}Z_{_{N(t+\tau)}}^{2}||\mathbf{e}||_{P,r}
+&\Big(\sum_{i=0}^{N(t+\tau)- N(t)-2}\big(1+Z_{_{N(t+\tau)}}^{2}||\mathbf{e}\mathbf{a}^{\prime}||_{P,r}\big)\times\\
&\quad\cdots\times\big(1+Z_{_{N(t+\tau)-i}}^{2}||\mathbf{e}\mathbf{a}^{\prime}||_{P,r}\big)
\times\alpha_{0}Z_{_{N(t+\tau)-i-1}}^{2}||\mathbf{e}||_{P,r}\Big).
\end{align*}
Since $E\big(Z_{i}^{2}\big)<\infty$ and $E\big(N(t+\tau)-N(t)\big)=\Lambda(t+\tau)-\Lambda(t)<\infty,$ it follows that $E\big(||\mathbf{K}_{t,t+\tau}||_{P,r}\big)<\infty$. So, by Jensen's inequality, we have that 
\begin{align*}
E\Big(log^{+}||\mathbf{K}_{t,t+\tau}||_{P,r}\Big)=E\Big(log\big(max\lbrace1,||\mathbf{K}_{t,t+\tau}||_{P,r}\rbrace\big)\Big)<\infty.
\end{align*}
Now, since $E\Big(log||J_{t,t+\tau}||_{P,r}\Big)<0$ and $E\Big(log^{+}||\mathbf{K}_{t,t+\tau}||_{P,r}\Big)<\infty$, it follows from the general theory of random recurrence equations (see \cite{b3}) that $H_{m}^{(t)}$, in (\ref{Th3.2(5.14)}), converges a.s. to 0 as $m\rightarrow\infty$ and that $\mathbf{U}_{m}^{(t)}$ converges a.s. absolutely to some random vector $\mathbf{U}^{(t)}$ as $m\rightarrow\infty$, for fixed $t\in[0,\tau)$. Using the fact that the process $(\mathbf{Y}_{t})_{t\geq0}$ has c\`adl\`ag paths, it follows that $||\mathbf{Y}_{t}||_{P,r}$ is a.s. finite for fixed $t\in[0,\tau)$. So
\begin{align*}
\lim_{m\rightarrow\infty}||H_{m}^{(t)}\mathbf{Y}_{t}||_{P,r}=0\qquad\textnormal{a.s.,}
\end{align*} 
and it follows the $\mathbf{Y}_{t+m\tau}$ converges in distribution to $\mathbf{U}^{(t)}$ as $m\rightarrow\infty$, for fixed $t\in[0,\tau)$.\\
It remains to show that $\mathbf{U}^{(t)}$ satisfies (\ref{equ3.3}) for fixed $t\in[0,\tau)$. It follows from (\ref{Th3.2(5.12)}) and the part (c) of Theorem 3.1 that $\mathbf{Y}_{t+(m-1)\tau}$ is independent of $(J_{t+(m-1)\tau,t+m\tau}$, $\mathbf{K}_{t+(m-1)\tau,t+m\tau})$. So, by the fact that $(J_{t+(m-1)\tau,t+m\tau},\mathbf{K}_{t+(m-1)\tau,t+m\tau})$ and $(J_{t,t+\tau},\mathbf{K}_{t,t+\tau})$ are i.i.d., for any $m\geq2$, and that $\mathbf{Y}_{t+m\tau}$ converges in distribution to $\mathbf{U}^{(t)}$ as $m\rightarrow\infty$, we have that
\begin{align}\label{Th3.2(5.22)}
\big(J_{t+(m-1)\tau,t+m\tau},\mathbf{K}_{t+(m-1)\tau,t+m\tau},\mathbf{Y}_{t+(m-1)\tau}\big)\overset{d}{\longrightarrow}
\big(J_{t,t+\tau},\mathbf{K}_{t,t+\tau},\mathbf{U}^{(t)}\big)\quad \textnormal{as}\quad m\rightarrow\infty,
\end{align}
where $\mathbf{U}^{(t)}$ is independent of $(J_{t,t+\tau},\mathbf{K}_{t,t+\tau})$ and $\overset{d}{\longrightarrow}$ denotes the convergence in distribution. Hence, by (\ref{Th3.2(5.12)}) and (\ref{Th3.2(5.22)}), it follows (\ref{equ3.3}).
\vspace{7mm}
\\
\textbf{P4: Proof of Corollary 3.3}
\vspace{2mm}

\textbf{Proof ($a$):} By Theorem 3.2, $\mathbf{U}^{(t)}$ is the long run behavior of $\mathbf{Y}_{t+m\tau}$ for fixed $t\in[0,\tau)$ and large $m$, so that $\mathbf{U}^{(t)}\overset{d}{=}J_{t,t+\tau}\mathbf{U}^{(t)}+\mathbf{K}_{t,t+\tau}$ and $\mathbf{U}^{(t)}$ is independent of $(J_{t,t+\tau},\mathbf{K}_{t,t+\tau})$. Since, by the part (b) of Theorem 3.1, $(J_{t,t+\tau}, \mathbf{K}_{t,t+\tau})$ and $(J_{t+(k-1)\tau,t+k\tau},\mathbf{K}_{t+(k-1)\tau,t+k\tau})$, for $k\geq2$, are i.i.d., we have that 
\begin{align}\nonumber
\mathbf{U}^{(t)}&\overset{d}{=}J_{t,t+\tau}\mathbf{U}^{(t)}+\mathbf{K}_{t,t+\tau}\\ \label{Coro3.3(5.23)}
&\overset{d}{=}J_{t+(k-1)\tau,t+k\tau}\mathbf{U}^{(t)}+\mathbf{K}_{t+(k-1)\tau,t+k\tau}.
\end{align}
It also follows from the parts ($a$) and (c) of Theorem 3.1 that, for fixed $t\in[0,\tau)$ and $k\in\mathbb{N}$,
\begin{align}\label{Coro3.3(5.24)}
\mathbf{Y}_{t+k\tau}=J_{t+(k-1)\tau,t+k\tau}\mathbf{Y}_{t+(k-1)\tau}+\mathbf{K}_{t+(k-1)\tau,t+k\tau}
\end{align}
and that $\mathbf{Y}_{t+(k-1)\tau}$ is independent of $(J_{t+(k-1)\tau,t+k\tau},\mathbf{K}_{t+(k-1)\tau,t+k\tau})$. Thus if $\mathbf{Y}_{t}\overset{d}{=}\mathbf{U}^{(t)}$ for $t\in[0,\tau)$, then (\ref{Coro3.3(5.23)}) and (\ref{Coro3.3(5.24)}) imply that 
$\mathbf{Y}_{t+\tau}\overset{d}{=}\mathbf{U}^{(t)}$. Further, it follows by induction on $k\in\mathbb{N}$ that
\begin{align}\label{Coro3.3(5.25)}
\mathbf{Y}_{t+k\tau}\overset{d}{=}\mathbf{U}^{(t)}.
\end{align}
We show that the process $(\mathbf{Y}_{t})_{t\geq0}$ is strictly periodically stationary with period $\tau$. If $n=1$, then 
for any $t_{_{1}}\geq0$ there exist $t\in[0,\tau)$ and $k\in\mathbb{N}^{0}$ such that $t_{_{1}}=t+k\tau$. So, by (\ref{Coro3.3(5.25)}) 
we have that
\begin{align*}
\mathbf{Y}_{t_{_{1}}}\overset{d}{=}\mathbf{Y}_{t_{_{1}}+\tau}.
\end{align*}
Now, let $n\geq2$ and $0\leq t_{_{1}}\leq t_{_{2}}\leq\cdots\leq t_{_{n}}$. Since, by the part ($a$) of Theorem 3.1,
\begin{align}\label{Coro3.3(5.26)}
\mathbf{Y}_{t_{_{n}}}=J_{t_{_{n-1}},t_{_{n}}}\mathbf{Y}_{t_{_{n-1}}}+\mathbf{K}_{t_{_{n-1}},t_{_{n}}},
\end{align}
it follows that $(\mathbf{Y}_{t})_{t\geq0}$ is a Markov process. So, by the fact that $(J_{t_{_{n-1}},t_{_{n}}},\mathbf{K}_{t_{_{n-1}},t_{_{n}}})\overset{d}{=}(J_{t_{_{n-1}}+\tau,t_{_{n}}+\tau}, \mathbf{K}_{t_{_{n-1}}+\tau,t_{_{n}}+\tau})$, we have that
\begin{align}\nonumber
\big(\mathbf{Y}_{t_{_{n}}}|\mathbf{Y}_{t_{_{n-1}}}=\mathbf{y}_{n-1},\cdots, \mathbf{Y}_{t_{_{1}}}=\mathbf{y}_{1}\big)
&=J_{t_{_{n-1}},t_{_{n}}}\mathbf{y}_{n-1}+\mathbf{K}_{t_{_{n-1}},t_{_{n}}}\\ \nonumber
&\overset{d}{=}J_{t_{_{n-1}}+\tau,t_{_{n}}+\tau}\mathbf{y}_{n-1}+\mathbf{K}_{t_{_{n-1}}+\tau,t_{_{n}}+\tau}\\ \label{Coro3.3(5.26.2)}
&=\big(\mathbf{Y}_{t_{_{n}}+\tau}|\mathbf{Y}_{t_{_{n-1}}+\tau}=\mathbf{y}_{n-1},\cdots, 
\mathbf{Y}_{t_{_{1}}+\tau}=\mathbf{y}_{1}\big).
\end{align}
Therefore, an induction argument for $n=2, 3, \cdots$ and (\ref{Coro3.3(5.26.2)}) imply
\begin{align*}
\big(\mathbf{Y}_{t_{_{1}}}, \mathbf{Y}_{t_{_{2}}}, \cdots, \mathbf{Y}_{t_{_{n}}}\big)\overset{d}
{=}\big(\mathbf{Y}_{t_{_{1}}+\tau}, \mathbf{Y}_{t_{_{2}}+\tau}, \cdots, \mathbf{Y}_{t_{_{n}}+\tau}\big).
\end{align*}
Hence, by (\ref{equ2.6}) it follows that
\begin{align}\label{Coro3.3(5.27)}
\big(V_{t_{_{1}}}, V_{t_{_{2}}}, \cdots, V_{t_{_{n}}}\big)\overset{d}
{=}\big(V_{t_{_{1}}+\tau}, V_{t_{_{2}}+\tau}, \cdots, V_{t_{_{n}}+\tau}\big).
\end{align}

\textbf{Proof (b):} By (\ref{equ2.6}) and (\ref{Th3.1(5.11)}) we have that
\begin{align}\nonumber
V_{u}&=\alpha_{0}+\mathbf{a}^{\prime}\mathbf{Y}_{u^{-}}\\ \label{Coro3.3(5.28)}
&=\alpha_{0}+\mathbf{a}^{\prime}J_{0,u^{-}}\mathbf{Y}_{0}+\mathbf{a}^{\prime}\mathbf{K}_{0,u^{-}}\quad \forall u\geq0.
\end{align}
As $(J_{0,u^{-}},\mathbf{K}_{0,u^{-}})$ are constructed from segments of the driving process $S_{t}$ in the interval $(0,u^{-}]$, see Theorem 3.1 (c), and that $\mathbf{Y}_{0}$ is independent of $(S_{t})_{t\geq0}$, assumed by (\ref{equ2.7}), it follows that $V_{u}$ and $dS_{u}=S_{u+du}-S_{u}$ are independent. So, (\ref{Coro3.3(5.27)}) and semi-L\'evy structure of $S_{t}$ imply that $V_{u}$ and $dS_{u}=S_{u+du}-S_{u}$ both are periodic with the same period $\tau$, say 
\begin{align}\label{Coro3.3(5.29)}
\sqrt{V_{u}}dS_{u}\overset{d}{=}\sqrt{V_{u+\tau}}dS_{u+\tau}.
\end{align}
Thus it follows from (\ref{Coro3.3(5.29)}) that
\begin{align}\nonumber
E\big(G_{t}^{(l)}\big)&=E\Big(\int_{t}^{t+l}\sqrt{V_{u}}dS_{u}\Big)\\ \label{Coro3.3(5.30)}
&=E\Big(\int_{t}^{t+l}\sqrt{V_{u+\tau}}dS_{u+\tau}\Big)=E\big(G_{t+\tau}^{(l)}\big)
\end{align}
and
\begin{align}\nonumber
E\big(G_{t}^{(l)}G_{t+h}^{(l)}\big)&=
E\Big(\int_{t}^{t+l}\sqrt{V_{u}}dS_{u}\int_{t+h}^{t+h+l}\sqrt{V_{u}}dS_{u}\Big)\\ \label{Coro3.3(5.31)}
&=E\Big(\int_{t}^{t+l}\sqrt{V_{u+\tau}}dS_{u+\tau}\int_{t+h}^{t+h+l}\sqrt{V_{u+\tau}}dS_{u+\tau}\Big)
=E\big(G_{t+\tau}^{(l)}G_{t+h+\tau}^{(l)}\big).
\end{align}
Hence, from (\ref{Coro3.3(5.30)}) and (\ref{Coro3.3(5.31)}) we have that
\begin{align*}
cov\big(G_{t}^{(l)},G_{t+h}^{(l)}\big)=cov\big(G_{t+\tau}^{(l)},G_{t+h+\tau}^{(l)}\big).
\end{align*}
\\
\textbf{P5: Proof of Corollary 3.4}
\vspace{2mm}

\textbf{Proof ($a$):} By the part (b) of Corollary 3.3, it follows that $V_{u}$ and $dS_{u}=S_{u+du}-S_{u}$ are independent for any $u\geq0$. So, $E(dS_{u})=E(S_{u+du})-E(S_{u})=0$ and 
\begin{align*}\nonumber
E\big(G_{t}^{(l)}\big)&=E\Big(\int_{t}^{t+l}\sqrt{V_{u}}dS_{u}\Big)\\
&=\int_{t}^{t+l}E(\sqrt{V_{u}})E(dS_{u})=0.
\end{align*}
Since, by Definition 2.2, the jump sizes $Z_{n}$ have finite second order moment it is easy to see that $E(S_{t})^{2}<\infty$. So, it follows from the proposition 3.17 of Cont and Tankov \cite{c} that $S_{t}$ is a martingale. Thus, from It\^o isometry for square integrable martingales as integrators (e.g. Rogers and Williams \cite{r1}, IV 27) follows
\begin{align*}\nonumber
E\big(G_{t}^{(l)}G_{t+h}^{(l)}\big)=E\Big(\int_{0}^{t+h+l}V_{u}I_{(t,t+l]}(u)I_{(t+h,t+h+l]}(u)d[S,S]_{u}\Big)=0
\end{align*}
for $h\geq l$. Hence, $cov\big(G_{t}^{(l)},G_{t+h}^{(l)}\big)=0$.
\vspace{3mm}

\textbf{Proof (b):} By (\ref{Coro3.3(5.29)}) we have that
\begin{align}\nonumber
E\Big((G_{t}^{(l)})^{2}\Big)&=E\Big(\int_{t}^{t+l}\sqrt{V_{u}}dS_{u}\Big)^{2}\\ \label{Coro3.4(5.32)}
&=E\Big(\int_{t}^{t+l}\sqrt{V_{u+\tau}}dS_{u+\tau}\Big)^{2}=E\Big((G_{t+\tau}^{(l)})^{2}\Big)
\end{align}
and
\begin{align}\nonumber
E\Big((G_{t}^{(l)})^{2}(G_{t+h}^{(l)})^{2}\Big)&=
E\Big(\big(\int_{t}^{t+l}\sqrt{V_{u}}dS_{u}\big)^{2}\big(\int_{t+h}^{t+h+l}\sqrt{V_{u}}dS_{u}\big)^{2}\Big)\\ \nonumber
&=E\Big(\big(\int_{t}^{t+l}\sqrt{V_{u+\tau}}dS_{u+\tau}\big)^{2}\big(\int_{t+h}^{t+h+l}\sqrt{V_{u+\tau}}dS_{u+\tau}\big)^{2}\Big)\\ 
\label{Coro3.4(5.33)}
&=E\Big((G_{t+\tau}^{(l)})^{2}(G_{t+h+\tau}^{(l)})^{2}\Big).
\end{align}
Hence, from (\ref{Coro3.4(5.32)}) and (\ref{Coro3.4(5.33)}) it follows that
\begin{align*}
cov\Big((G_{t}^{(l)})^{2},(G_{t+h}^{(l)})^{2}\Big)=cov\Big((G_{t+\tau}^{(l)})^{2},(G_{t+h+\tau}^{(l)})^{2}\Big).
\end{align*}
\\
\textbf{P6: Proof of Theorem 3.5}\\
Suppose that (\ref{equ3.4}) and (\ref{equ3.5}) both hold. The following proof of Theorem 3.1, for any $t\geq0$ it follows that $t\in[\Upsilon_{_{N(t)}},\Upsilon_{_{N(t)+1}})$. So, we have by (\ref{Th3.1(5.6)}) that 
\begin{align}\label{Th3.5(5.34)}
\mathbf{Y}_{t}=e^{B(t-\Upsilon_{_{N(t)}})}\mathbf{Y}_{\Upsilon_{_{N(t)}}}.
\end{align}
By (\ref{equ2.6}) and replacing from $N(t)$ in $n$, (\ref{Th3.1(5.3)}) reads as 
\begin{align}\label{Th3.5(5.35)}
\mathbf{Y}_{\Upsilon_{_{N(t)}}}=\mathbf{Y}_{\Upsilon^{-}_{_{N(t)}}}+\mathbf{e}V_{\Upsilon_{_{N(t)}}}Z_{_{N(t)}}^{2}
\end{align}
and (\ref{Th3.1(5.4)}) as
\begin{align}\label{Th3.5(5.36)}
\mathbf{Y}_{\Upsilon^{-}_{_{N(t)}}}=e^{B(\Upsilon_{_{N(t)}}-\Upsilon_{_{N(t)-1}})}\mathbf{Y}_{\Upsilon_{_{N(t)-1}}}.
\end{align}
Replacing (\ref{Th3.5(5.36)}) in (\ref{Th3.5(5.35)}) implies that
\begin{align}\label{Th3.5(5.37)}
\mathbf{Y}_{\Upsilon_{_{N(t)}}}=e^{B(\Upsilon_{_{N(t)}}-\Upsilon_{_{N(t)-1}})}\mathbf{Y}_{\Upsilon_{_{N(t)-1}}}
+\mathbf{e}V_{\Upsilon_{_{N(t)}}}Z_{_{N(t)}}^{2}.
\end{align}
So, it follows from (\ref{Th3.5(5.34)}) and (\ref{Th3.5(5.37)}) that
\begin{align}\nonumber
\mathbf{Y}_{t}&=e^{B(t-\Upsilon_{_{N(t)}})}\Big[e^{B(\Upsilon_{_{N(t)}}-\Upsilon_{_{N(t)-1}})}\mathbf{Y}_{\Upsilon_{_{N(t)-1}}}
+\mathbf{e}V_{\Upsilon_{_{N(t)}}}Z_{_{N(t)}}^{2}\Big]\\ \label{Th3.5(5.38)}
&=e^{B(t-\Upsilon_{_{N(t)-1}})}\mathbf{Y}_{\Upsilon_{_{N(t)-1}}}
+e^{B(t-\Upsilon_{_{N(t)}})}\mathbf{e}V_{\Upsilon_{_{N(t)}}}Z_{_{N(t)}}^{2}.
\end{align}
By successive replacement from (\ref{Th3.5(5.37)}) in (\ref{Th3.5(5.38)}) we have that
\begin{align*}
\mathbf{Y}_{t}=e^{Bt}\mathbf{Y}_{0}+\sum_{i=1}^{N(t)}e^{B(t-\Upsilon_{_{i}})}\mathbf{e}V_{\Upsilon_{_{i}}}Z_{_{i}}^{2},
\qquad t\geq0.
\end{align*}
So, from this and (\ref{equ3.5}) follows that 
\begin{align}\label{Th3.5(5.39)}
\mathbf{a}^{\prime}\mathbf{Y}_{t}&=\mathbf{a}^{\prime} e^{Bt}\mathbf{Y}_{0}
+\sum_{i=1}^{N(t)}\mathbf{a}^{\prime}e^{B(t-\Upsilon_{_{i}})}\mathbf{e}V_{\Upsilon_{_{i}}}Z_{_{i}}^{2}\\ \label{Th3.5(5.40)}
&\geq\gamma+\sum_{i=1}^{N(t)}\mathbf{a}^{\prime}e^{B(t-\Upsilon_{_{i}})}\mathbf{e}V_{\Upsilon_{_{i}}}Z_{_{i}}^{2}.
\end{align}
By a similar method in Brockwell et al. \cite{b5}, we find by (\ref{equ2.6}) and (\ref{Th3.5(5.40)}) that $V_{t}=\alpha_{0}+\mathbf{a}^{\prime}\mathbf{Y}_{t^{-}}\geq\alpha_{0}+\gamma$ for $t\in[0,\Upsilon_{_{1}}]$. It also follows from (\ref{Th3.5(5.40)}), (\ref{equ3.5}) and (\ref{equ3.4}) that
\begin{align*}
V_{\Upsilon^{+}_{_{1}}}=\alpha_{0}+\mathbf{a}^{\prime}\mathbf{Y}_{\Upsilon_{_{1}}}
&=\alpha_{0}+\mathbf{a}^{\prime}e^{B\Upsilon_{_{1}}}\mathbf{Y}_{0}+
\mathbf{a}^{\prime}\mathbf{e}V_{\Upsilon_{_{1}}}Z_{_{1}}^{2}\\
&\geq\alpha_{0}+\gamma.
\end{align*}
So, by induction one can easily verify $V_{t}\geq\alpha_{0}+\gamma\geq0$ for any $t\geq0$.
\vspace{2mm}

To prove the converse, suppose first that (\ref{equ3.5}) fails. Then, it follows from (\ref{Th3.5(5.39)}) that, for  $t\in[0,\Upsilon_{_{1}}]$,
\begin{align}\label{Th3.5(5.41)}
V_{t}=\alpha_{0}+\mathbf{a}^{\prime}\mathbf{Y}_{t^{-}}<\alpha_{0}+\gamma.
\end{align}
Since $P(\Upsilon_{_{1}}>0)>0$, it follows from (\ref{Th3.5(5.41)}) and $\gamma=-\alpha_{0}$ that $P(V_{t}<0)>0$ for $t\in[0,\Upsilon_{_{1}}]$. Now assume that (\ref{equ3.5}) holds with $\gamma>-\alpha_{0}$, but (\ref{equ3.4}) fails. Suppose that the support of the jumps distribution $Z_{i}$ is unbounded. Let $(t_{1},t_{2})\subset(0,\infty)$ be an interval such that $\mathbf{a}^{\prime}e^{Bt}\mathbf{e}<c$ for any $t\in(t_{1},t_{2})$ and for some $c<0$. By (\ref{equ2.6}), (\ref{Th3.5(5.39)}) and (\ref{equ3.5}) we have that
\begin{align}\nonumber
V_{\Upsilon_{_{1}}}=\alpha_{0}+\mathbf{a}^{\prime}\mathbf{Y}_{\Upsilon_{_{1}}^{-}}&=\alpha_{0}
+\mathbf{a}^{\prime}e^{B\Upsilon_{_{1}}}\mathbf{Y}_{0}\\ 
\label{Th3.5(5.42)}
&\geq\alpha_{0}+\gamma.
\end{align}
It follows from (\ref{Th3.5(5.42)}) and assumption $\gamma>-\alpha_{0}$, that $P(V_{\Upsilon_{_{1}}}\geq\alpha_{0}+\gamma>0)=1$. So, one can easily verify that the set 
\begin{align*}\nonumber
A:=\lbrace \Upsilon_{_{1}}<t_{3}<\Upsilon_{_{2}}, t_{3}-\Upsilon_{_{1}}\in(t_{1}, t_{2}), V_{\Upsilon_{_{1}}}>0\rbrace
\end{align*}
has positive probability for $t_{3}>t_{2}$. On $A$, by (\ref{Th3.5(5.39)}) we have that
\begin{align}\label{Th3.5(5.43)}
V_{t_{3}}=\alpha_{0}+\mathbf{a}^{\prime}e^{Bt_{3}}\mathbf{Y}_{0}+\mathbf{a}^{\prime}e^{B(t_{3}-\Upsilon_{_{1}})}\mathbf{e}V_{\Upsilon_{_{1}}}Z_{1}^{2}.
\end{align}
Since $\mathbf{a}^{\prime}e^{B(t_{3}-\Upsilon_{_{1}})}\mathbf{e}<c<0$, by choosing $Z_{1}$ large enough it follows from (\ref{Th3.5(5.43)}) that $P(V_{t_{3}}<0)>0$.

\end{document}